\RequirePackage{fix-cm}
\documentclass[smallcondensed]{svjour3}     
\smartqed  

\usepackage{mathptmx}      
\usepackage{graphicx}
\usepackage{amsmath}
\usepackage{amssymb}
\usepackage{setspace}
\usepackage{enumerate}
\usepackage[ruled,linesnumbered,vlined]{algorithm2e}
\graphicspath{{figs/}}
\usepackage{color}
\usepackage{array}
\newcolumntype{C}[1]{>{\centering\let\newline\\\arraybackslash\hspace{0pt}}m{#1}}
\usepackage{nicefrac}
\usepackage{microtype}
\usepackage{hyperref}
\usepackage{bm}

\definecolor{webgreen}{rgb}{0,.35,0}
\definecolor{webbrown}{rgb}{.6,0,0}
\definecolor{RoyalBlue}{rgb}{0,0,0.9}
\definecolor{purp}{rgb}{0.6,0.05,0.8}
\definecolor{ora}{rgb}{0.7,0.35,0.02}

\hypersetup{
   colorlinks=true, linktocpage=true, pdfstartpage=3, pdfstartview=FitV,
   breaklinks=true, pdfpagemode=UseNone, pageanchor=true, pdfpagemode=UseOutlines,
   plainpages=false, bookmarksnumbered, bookmarksopen=true, bookmarksopenlevel=1,
   hypertexnames=true, pdfhighlight=/O,
   urlcolor=webbrown, linkcolor=RoyalBlue, citecolor=webgreen,
   pdfauthor={Gary P. T. Choi},
   pdfsubject={Volumetric density-equalizing reference map with applications}
}

\newcounter{exctr}

\newcommand{\p}{\partial}
\newcommand{\R}{\mathbb{R}}
\newcommand{\bX}{\mathbf{X}}
\newcommand{\bx}{\mathbf{x}}
\newcommand{\bxi}{\boldsymbol\xi}
\newcommand{\bchi}{\boldsymbol\chi}

\journalname{}

\begin{document}

\title{Volumetric density-equalizing reference map with applications}


\author{Gary P. T. Choi \and Chris H. Rycroft}

\institute{
Gary P. T. Choi \at
John A. Paulson School of Engineering and Applied Sciences, Harvard University, Cambridge, MA 02138, USA \\
\email{pchoi@g.harvard.edu}           \\
\href{https://orcid.org/0000-0001-5407-9111}{\includegraphics[scale=0.45]{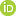}} \href{https://orcid.org/0000-0001-5407-9111}{https://orcid.org/0000-0001-5407-9111}
\and
Chris H. Rycroft \at
John A. Paulson School of Engineering and Applied Sciences, Harvard University, Cambridge, MA 02138, USA, and\\
Mathematics Group, Lawrence Berkeley National Laboratory, Berkeley, CA 94720, USA\\
\email{chr@seas.harvard.edu}           \\
\href{https://orcid.org/0000-0003-4677-6990}{\includegraphics[scale=0.45]{orcid}} \href{https://orcid.org/0000-0003-4677-6990}{https://orcid.org/0000-0003-4677-6990}
}

\date{Received: date / Accepted: date}

\maketitle
\begin{abstract}
The density-equalizing map, a technique developed for cartogram creation, has been widely applied to data visualization but only for 2D applications. In this work, we propose a novel method called the volumetric density-equalizing reference map (VDERM) for computing density-equalizing map for volumetric domains. Given a prescribed density distribution in a volumetric domain in $\mathbb{R}^3$, the proposed method continuously deforms the domain, with different volume elements enlarged or shrunk according to the density distribution. With the aid of the proposed method, medical and sociological data can be visualized via deformations of 3D objects. The method can also be applied to adaptive remeshing and shape modeling. Furthermore, by exploiting the time-dependent nature of the proposed method, applications to shape morphing can be easily achieved. Experimental results are presented to demonstrate the effectiveness of the proposed method.


\end{abstract}

\section{Introduction}
In the notable work by Gastner and Newman~\cite{Gastner04}, an algorithm for producing density-equalizing map projections was proposed based on the diffusion equation. Specifically, given a planar geographical map and a density distribution prescribed on every part of the map, the algorithm continuously deforms the map such that the difference in the density at different regions is transformed into a difference in the area of the regions. Regions with a larger prescribed density expand and those with a smaller density shrink. Ultimately, the density is equalized over the entire deformed map. The algorithm has been extensively applied to the visualization of biological and sociological data, such as the global amphibian species diversity~\cite{Wake08}, the global population and income~\cite{Dorling10}, the world citation network~\cite{Pan12}, the national climate contributions to observed global warming~\cite{Matthews14}, the global burden of drug-resistant tuberculosis in children~\cite{Dodd16}, and the impact of austerity and the economic crisis in Europe~\cite{Ballas17}.

In recent years, a few improvements and extensions of the above-mentioned work have been proposed. Choi and Rycroft~\cite{Choi18} developed a method for computing density-equalizing maps for simply-connected open surfaces and explored the close connection between map projections and surface parameterization, which opens up a wide range of applications of density-equalizing maps to geometry processing. Gastner et al.~\cite{Gastner18} proposed a new approach for computing density-equalizing map projections by considering a linearization of the density diffusion process. Recently, Choi et al.~\cite{Choi20} proposed a method for computing area-preserving density-equalizing maps for carotid artery flattening.

Note that all the above-mentioned works only focused on the formulation and application of density-equalizing map in 2D. With the advancement of technology, there is an increasing need for 3D data processing. In this work, we develop a novel method called the \emph{volumetric density-equalizing reference map} (VDERM) that produces volumetric deformations based on a prescribed density distribution. Our method combines and extends the 2D density-equalizing map~\cite{Gastner04} and the reference map technique in solid mechanics~\cite{Kamrin12,Valkov15,Rycroft18}. The proposed method can be utilized for volumetric data visualization, shape modeling, adaptive remeshing etc. To the best of our knowledge, this is the first work on higher dimensional density-equalizing maps. Furthermore, all the previous works only focused on the use of density-equalizing map for producing a single final output. With the observation that density-equalizing map is a continuous deformation, we introduce the use of density-equalizing map for time-dependent applications such as shape morphing.

The rest of the paper is organized as follows. In Section \ref{sect:contributions}, we highlight the contributions of our work. In Section \ref{sect:background}, we review the background of diffusion-based map-making and the reference map technique in 2D. In Section \ref{sect:volumetric}, we describe our proposed method for volumetric density-equalizing maps and present numerical experiments to demonstrate its effectiveness. We then discuss several novel applications of volumetric density-equalizing reference map in Section \ref{sect:application}. In Section \ref{sect:discussion}, we conclude our work and discuss possible future works.

\section{Contributions} \label{sect:contributions}
The contributions of our work are three-fold:
\begin{enumerate}[(i)]
 \item Previous works of density-equalizing maps have only considered the formulation and its applications in 2D. Our work is the first work focusing on higher dimensional density-equalizing maps. Experimental results show that our method is capable of producing volumetric deformations accurately based on the prescribed density. 
 \item Our proposed volumetric method leads to novel applications of density-equalizing map, such as 3D medical and sociological data visualization, adaptive volumetric remeshing, and deformation-based shape modeling.
 \item All the previous works on density-equalizing maps have only focused on the use of the final result produced by the density-equalizing map method. Our work is the first work that considers the use of the entire density-equalization process, i.e. not only the final maps but also the intermediate states obtained throughout the process, for time-dependent applications such as object morphing. 
\end{enumerate}

\section{Background} \label{sect:background}
\subsection{Diffusion-based cartogram}
\begin{figure}[t]
\centering
\includegraphics[width=0.8\textwidth]{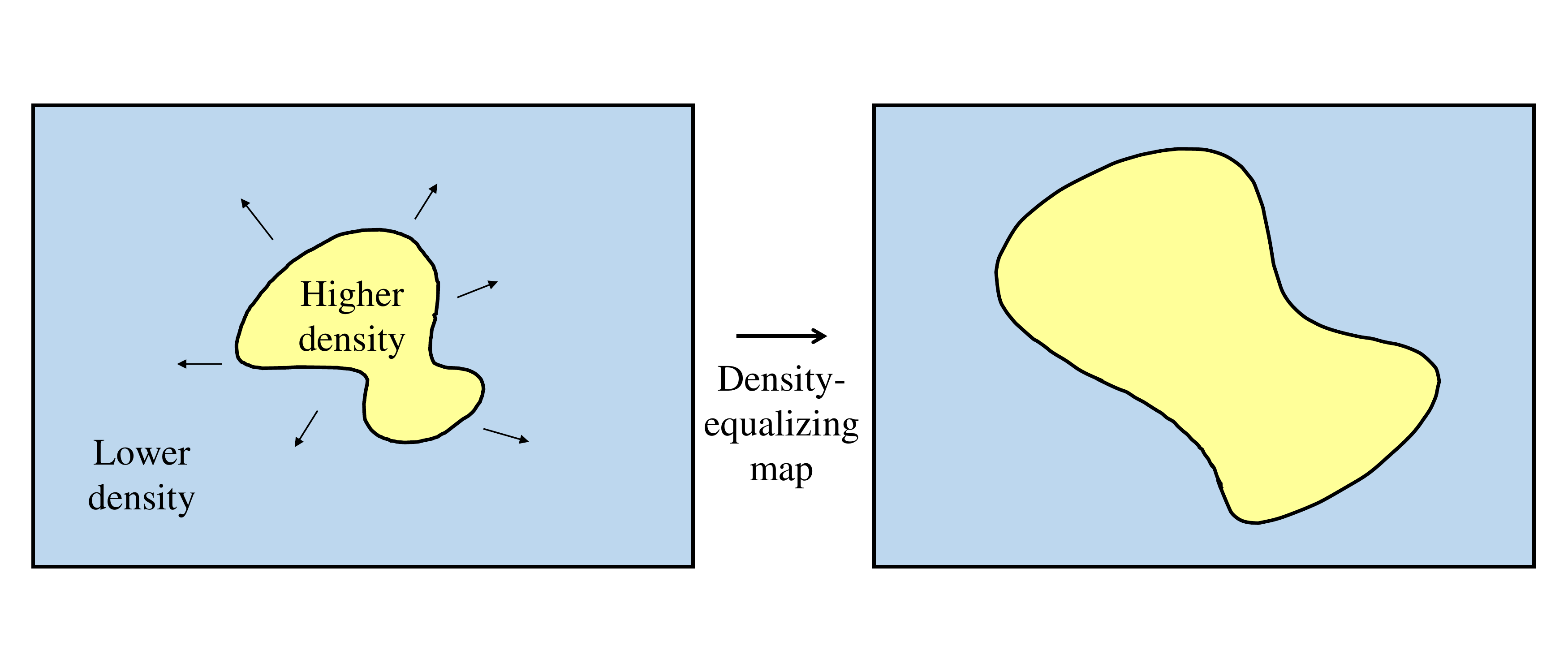}
\caption{An illustration of the diffusion-based cartogram. The density diffusion creates a velocity field that enlarges regions with higher density and shrinks regions with lower density.}
\label{fig:dem_illustration}
\end{figure}

Gastner and Newman~\cite{Gastner04} proposed a method for creating cartograms based on density diffusion. Given a 2D domain, a positive quantity $\rho$ called the \emph{density} is first prescribed on every part of it. For instance, the 2D domain can be a geographical map of a country, and $\rho$ can be the population density of different provinces. The goal of a diffusion-based cartogram is to deform the domain by enlarging the regions with higher $\rho$, and shrinking the regions with lower $\rho$. This is achieved by equalizing $\rho$, following the advection equation
\begin{equation}
\frac{\partial \rho}{\partial t} = -\nabla \cdot \mathbf{j}, 
\end{equation}
where $\mathbf{j} = - \nabla \rho$ is the density flux. This yields the diffusion equation
\begin{equation}\label{eqt:diffusion}
\frac{\partial \rho}{\partial t} = \Delta \rho. 
\end{equation}
Any tracer particles on the domain will then move with velocity
\begin{equation}\label{eqt:velocity_field}
\mathbf{v} = - \frac{\nabla \rho}{\rho}.
\end{equation}
Therefore, the location of a tracer particle $\mathbf{x}$ at time $t$ can be given by
\begin{equation}\label{eqt:displacement}
\mathbf{x}(t) = \mathbf{x}(0) + \int_0^t \mathbf{v}(\mathbf{x},\tau) d \tau.
\end{equation}
As $t \to \infty$, $\rho$ is fully equalized on the 2D domain, and the resulting deformation produces an accurate density-based deformation according to the prescribed $\rho$. An illustration is given in Fig.~\ref{fig:dem_illustration}.

\begin{figure}[t]
\centering
\includegraphics[width=0.5\textwidth]{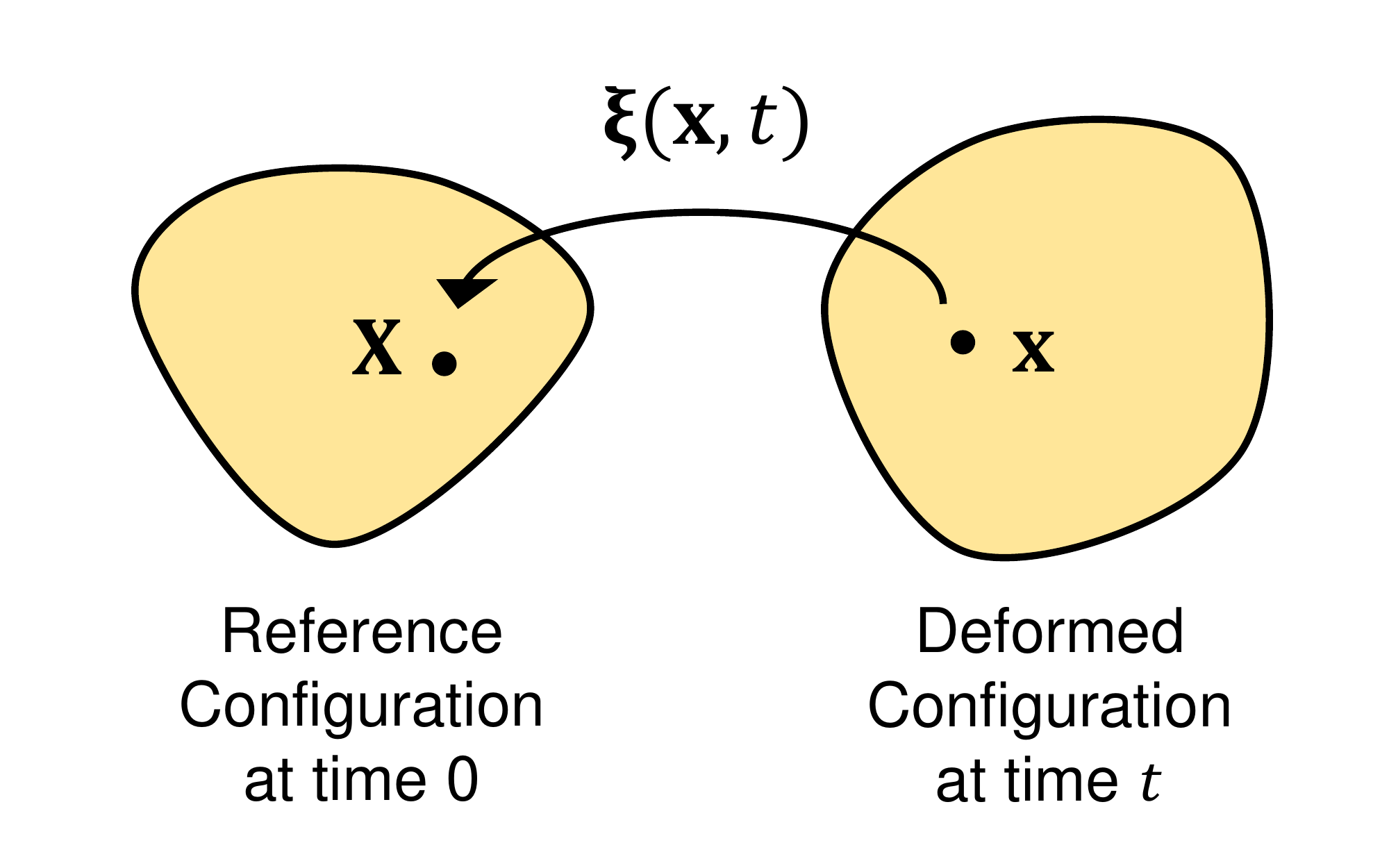}
\caption{An illustration of the reference map. ${\bm \xi}(\mathbf{x},t)$ is the inverse of the motion function, which indicates the reference location $\mathbf{X}$ of the material occupying the position $\mathbf{x}$ at time $t$.}
\label{fig:illustration_reference_map}
\end{figure}

\subsection{The reference map technique}
In solid mechanics, a standard mathematical approach is to consider a body in an undeformed reference configuration at time 0 that is mapped into a deformed configuration at time $t$~\cite{Gurtin10}. Let $\bX \in \R^d$ be a point in the reference configuration, and $\bx \in \R^d$ be the corresponding point in the deformed configuration. Then the deformation can be described by the \textit{motion function} $\bchi(\bX,t)=\bx$. The \textit{reference map} is then defined as the inverse motion function $\bxi(\bx,t)=\bX$ (Fig.~\ref{fig:illustration_reference_map}), which can be regarded as a vector field in the deformed configuration indicating the reference location of the material occupying the position $\bx$ at time $t$. It has been used in various applications such as for inverse design~\cite{Govindjee96,Fachinotti08}.

Kamrin, Rycroft, and coworkers developed the reference map technique (RMT)~\cite{Kamrin12,Valkov15,Rycroft18}, which uses the reference map as the basis for a fully Eulerian approach to solid mechanics, since it provides a simple and effective method to describe an arbitrary deformation of a body. Since the initial configuration is undeformed, we have $\bxi(\bx,0)=\bx$. Now, for any tracer particle, the reference location of it is the same at all time $t$ under the deformation. Therefore, we have
\begin{equation}
  \dot{\bxi} = \mathbf{0},
\end{equation}
yielding the advection equation
\begin{equation}\label{eqt:reference_map}
  \frac{\p \bxi}{\p t} + \mathbf{v} \cdot \nabla \bxi = \mathbf{0}.
\end{equation}
Eq.~\eqref{eqt:reference_map} can be solved on a discrete set of fixed Eulerian grid points using simple finite-difference schemes. In the RMT $\bxi$ is used to calculate mechanical quantities such as the Cauchy stress, but here we use it as a straightforward Eulerian-frame calculation of solid deformations.

\section{Volumetric density-equalizing reference map (VDERM)} \label{sect:volumetric}
Based on the 2D diffusion-based cartogram and the reference map technique introduced above, we propose a method for computing the volumetric density-equalizing reference map of 3D domains.
\subsection{Formulation}
Suppose we are given a rectangular solid domain $D \subset \mathbb{R}^3$, discretized as a $L\times M\times N$ 3D grid with grid spacing $h$ in all the $x$-, $y$- and $z$-directions. Let the coordinates of the grid points be $(ih, jh, kh)$, where $0 \leq i \leq L - 1$, $0 \leq j \leq M-1$, $0 \leq k \leq N-1$.  

\begin{figure}[t]
\centering
\includegraphics[width=\textwidth]{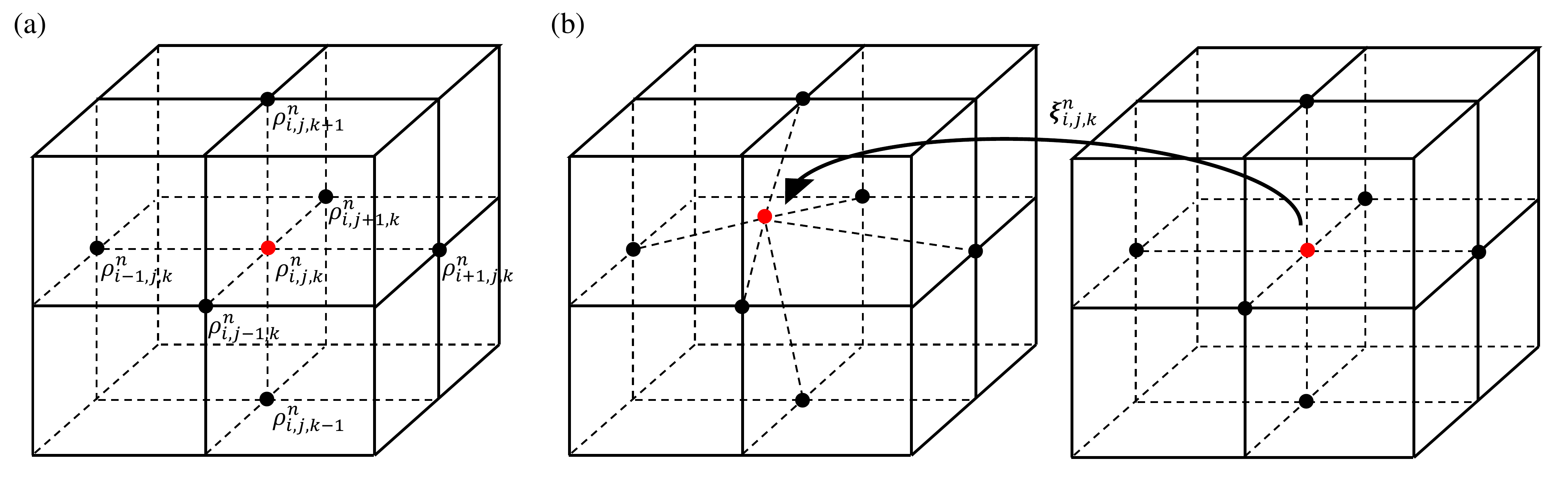}
\caption{A illustration of density diffusion and reference map update in our proposed volumetric density-equalizing reference map method. (a) The diffusion of $\rho$ is done iteratively by solving the diffusion equation~\eqref{eqt:derm_diffusion_discrete}, based on the six neighboring nodes. (b) The density gradient is then used for updating the reference map field $\bm \xi$ via the advection equation~\eqref{eqt:derm_advection_discrete}.}
\label{fig:vderm_illustration}
\end{figure}

Consider a prescribed density $\rho^0 = \rho^0(\mathbf{x})$ defined in $D$, discretized as $\rho_{i,j,k}^0 = \rho^0(ih, jh, kh)$ for all $i,j,k$. Let $\rho(\mathbf{x},t)$ be the density at time $t$, with $\rho(\mathbf{x},0) = \rho^0(\mathbf{x})$. The diffusion of $\rho$ follows from the diffusion equation
\begin{equation} \label{eqt:derm_diffusion}
\frac{\partial \rho}{\partial t}(\mathbf{x},t) = \Delta \rho (\mathbf{x},t),
\end{equation} 
which can be discretized using the following backward Euler scheme (Fig.~\ref{fig:vderm_illustration}):
\begin{equation} \label{eqt:derm_diffusion_discrete}
\begin{split}
\frac{\rho_{i,j,k}^{n} - \rho_{i,j,k}^{n-1}}{\delta t} 
=& \frac{1}{h^2} \left(\rho_{i+1,j,k}^n+\rho_{i-1,j,k}^n+\rho_{i,j+1,k}^n \right. + \left. \rho_{i,j-1,k}^n+\rho_{i,j,k+1}^n+\rho_{i,j,k-1}^n-6\rho_{i,j,k}^n \right),
\end{split}
\end{equation}
where $\delta t$ is the timestep to be determined, and $n = t / \delta t$ is the iteration number. Eq.~\eqref{eqt:derm_diffusion_discrete} can be further simplified as
\begin{equation} \label{eqt:derm_diffusion_discrete_simplified}
\rho^n = (I - \delta t \Delta)^{-1} \rho^{n-1},
\end{equation}
where $\Delta$ is an $LMN \times LMN$ matrix representing the seven-point Laplacian stencil. Note that $I - \delta t \Delta$ is a sparse symmetric positive definite matrix, and is only required to be precomputed once as it is independent of $n$.

Throughout the diffusion process, the density gradient of $\rho$ induces a velocity field 
\begin{equation} \label{eqt:derm_velocity}
\mathbf{v}(\mathbf{x},t) = - \frac{\nabla \rho (\mathbf{x},t)}{\rho (\mathbf{x},t)},
\end{equation}
which can be discretized using the central difference scheme:
\begin{equation} \label{eqt:derm_velocity_discrete}
\left\{\begin{split}
(\mathbf{v}_x)_{i,j,k}^n &= -\frac{\rho_{i+1,j,k}^{n}-\rho_{i-1,j,k}^{n}}{2h \rho_{i,j,k}^{n}},\\
(\mathbf{v}_y)_{i,j,k}^n &= -\frac{\rho_{i,j+1,k}^{n}-\rho_{i,j-1,k}^{n}}{2h \rho_{i,j,k}^{n}},\\
(\mathbf{v}_z)_{i,j,k}^n &= -\frac{\rho_{i,j,k+1}^{n}-\rho_{i,j,k-1}^{n}}{2h \rho_{i,j,k}^{n}}.
\end{split}\right.
\end{equation}

To track the deformation of $D$ under the velocity field, consider the reference map ${\bm \xi}(\mathbf{x},t)$, discretized as ${\bm \xi}_{i,j,k}^{n} = {\bm \xi}((ih, jh, kh), n \delta t)$. ${\bm \xi}$ can be obtained by solving the advection equation
\begin{equation} \label{eqt:derm_advection}
\frac{\p \bm \xi}{\p t}(\mathbf{x},t) - \frac{\nabla \rho (\mathbf{x},t)}{\rho (\mathbf{x},t)} \cdot \nabla \bm \xi(\mathbf{x},t) = \mathbf{0},
\end{equation}
which can be discretized using the second-order upwind scheme:
\begin{equation} \label{eqt:derm_advection_discrete}
\frac{{\bm \xi}_{i,j,k}^{n} - {\bm \xi}_{i,j,k}^{n-1}}{\delta t} = d_x + d_y + d_z,
\end{equation}
where 
\begin{equation}
d_x = \left\{\begin{array}{ll}
(\mathbf{v}_x)_{i,j,k}^{n} d_x^- & \text{ if } (\mathbf{v}_x)_{i,j,k}^{n}>0,\\
(\mathbf{v}_x)_{i,j,k}^{n} d_x^+ & \text{ if } (\mathbf{v}_x)_{i,j,k}^{n}\leq 0,
\end{array}\right.
\end{equation}
\begin{equation}
d_y = \left\{\begin{array}{ll}
(\mathbf{v}_y)_{i,j,k}^{n} d_y^- & \text{ if } (\mathbf{v}_y)_{i,j,k}^{n}>0,\\
(\mathbf{v}_y)_{i,j,k}^{n} d_y^+ & \text{ if } (\mathbf{v}_y)_{i,j,k}^{n}\leq 0,
\end{array}\right.
\end{equation}
\begin{equation}
d_z = \left\{\begin{array}{ll}
(\mathbf{v}_z)_{i,j,k}^{n} d_z^- & \text{ if } (\mathbf{v}_z)_{i,j,k}^{n}>0,\\
(\mathbf{v}_z)_{i,j,k}^{n} d_z^+ & \text{ if } (\mathbf{v}_z)_{i,j,k}^{n} \leq 0,
\end{array}\right.
\end{equation}
with $d_x^-, d_x^+, d_y^-, d_y^+, d_z^-, d_z^+$ being the three-point finite-difference discretization of the spatial derivative $\nabla \bm \xi$:
\begin{equation}
\left\{\begin{array}{ll}
d_x^- &= \frac{3{\bm \xi}_{i,j,k}^{n-1} - 4{\bm \xi}_{i-1,j,k}^{n-1} + {\bm \xi}_{i-2,j,k}^{n-1}}{2h},\\
d_x^+ &= \frac{-{\bm \xi}_{i+2,j,k}^{n-1} + 4{\bm \xi}_{i+1,j,k}^{n-1} - 3{\bm \xi}_{i,j,k}^{n-1}}{2h},\\
d_y^- &= \frac{3{\bm \xi}_{i,j,k}^{n-1} - 4{\bm \xi}_{i,j-1,k}^{n-1} + {\bm \xi}_{i,j-2,k}^{n-1}}{2h},\\
d_y^+ &= \frac{-{\bm \xi}_{i,j+2,k}^{n-1} + 4{\bm \xi}_{i,j+1,k}^{n-1} - 3{\bm \xi}_{i,j,k}^{n-1}}{2h},\\
d_z^- &= \frac{3{\bm \xi}_{i,j,k}^{n-1} - 4{\bm \xi}_{i,j,k-1}^{n-1} + {\bm \xi}_{i,j,k-2}^{n-1}}{2h},\\
d_z^+ &= \frac{-{\bm \xi}_{i,j,k+2}^{n-1} + 4{\bm \xi}_{i,j,k+1}^{n-1} -3{\bm \xi}_{i,j,k}^{n-1}}{2h}.
\end{array}\right.
\end{equation}
For $i=1, L-2$, we replace $d_x^-$ and $d_x^+$ with the two-point difference
\begin{equation} 
\left\{\begin{array}{ll}
d_x^- &= \frac{{\bm \xi}_{i,j,k}^{n-1} - {\bm \xi}_{i-1,j,k}^{n-1}}{h},\\
d_x^+ &= \frac{{\bm \xi}_{i+1,j,k}^{n-1} - {\bm \xi}_{i,j,k}^{n-1}}{h}.
\end{array}\right.
\end{equation}
Similarly, for $j = 1, M-2$, we use
\begin{equation} 
\left\{\begin{array}{ll}
d_y^- &= \frac{{\bm \xi}_{i,j,k}^{n-1} - {\bm \xi}_{i,j-1,k}^{n-1}}{h},\\
d_y^+ &= \frac{{\bm \xi}_{i,j+1,k}^{n-1} - {\bm \xi}_{i,j,k}^{n-1}}{h},
\end{array}\right.
\end{equation}
and for $k = 1, N-2$ we use
\begin{equation} 
\left\{\begin{array}{ll}
d_z^- &= \frac{{\bm \xi}_{i,j,k}^{n-1} - {\bm \xi}_{i,j,k-1}^{n-1}}{h},\\
d_z^+ &= \frac{{\bm \xi}_{i,j,k+1}^{n-1} - {\bm \xi}_{i,j,k}^{n-1}}{h}.
\end{array}\right.
\end{equation}

As for the choice of the timestep $\delta t$, note that the backward Euler scheme~\eqref{eqt:derm_diffusion_discrete} for the diffusion equation is unconditionally stable. A necessary condition for the convergence of the second-order upwind scheme~\eqref{eqt:derm_advection_discrete} for the advection equation is given by the Courant-Friedrichs-Lewy (CFL) condition~\cite{Courant28}, which says that the numerical domain of dependence must include the physical domain of dependence. More specifically, consider multiplying $\delta t$ on both sides in Eq.~\eqref{eqt:derm_advection_discrete} and making ${\bm \xi}_{i,j,k}^{n}$ as the subject. Then, the remaining terms involving ${\bm \xi}_{i,j,k}^{n-1}$ will be
\begin{equation}
\begin{split}
&{\bm \xi}_{i,j,k}^{n-1}- \delta t \left(\frac{3{\bm \xi}_{i,j,k}^{n-1} |(\mathbf{v}_x)_{i,j,k}^{n}|}{2h} + \frac{3{\bm \xi}_{i,j,k}^{n-1} |(\mathbf{v}_y)_{i,j,k}^{n}|}{2h} + \frac{3{\bm \xi}_{i,j,k}^{n-1} |(\mathbf{v}_z)_{i,j,k}^{n}|}{2h}\right) \\= &{\bm \xi}_{i,j,k}^{n-1} \left(1-\delta t \left(\frac{3 |(\mathbf{v}_x)_{i,j,k}^{n}|}{2h} + \frac{3 |(\mathbf{v}_y)_{i,j,k}^{n}|}{2h} + \frac{3 |(\mathbf{v}_z)_{i,j,k}^{n}|}{2h}\right) \right).
\end{split}
\end{equation}
Note that we have the absolute signs in the above expression because for a negative $(\mathbf{v}_x)_{i,j,k}^{n} (\mathbf{v}_y)_{i,j,k}^{n}$, or $(\mathbf{v}_z)_{i,j,k}^{n}$, the corresponding $d_x^+, d_y^+$, or $d_z^+$ involves a negative sign in the coefficient of ${\bm \xi}_{i,j,k}^{n-1}$. Now, since the full numerical domain of dependence of the scheme must contain the physical domain of dependence, for all $i,j,k,n$, the coefficient of ${\bm \xi}_{i,j,k}^{n-1}$ must satisfy
\begin{equation}
0 \leq 1-\delta t \left(\frac{3 |(\mathbf{v}_x)_{i,j,k}^{n}|}{2h} + \frac{3 |(\mathbf{v}_y)_{i,j,k}^{n}|}{2h} + \frac{3 |(\mathbf{v}_z)_{i,j,k}^{n}|}{2h}\right) \leq 1,
\end{equation}
which implies that
\begin{equation}
\delta t \left(\frac{3|(\mathbf{v}_x)_{i,j,k}^{n}|}{2h} + \frac{3|(\mathbf{v}_y)_{i,j,k}^{n}|}{2h} + \frac{3|(\mathbf{v}_z)_{i,j,k}^{n}|}{2h}\right) \leq 1.
\end{equation}
In our problem, the magnitude of the velocity field decreases as $n$ increases because diffusion smooths out sharp gradients in density. Therefore, we should have
\begin{equation}
\delta t \leq \frac{2h}{3 \max_{i,j,k}(|(\mathbf{v}_x)_{i,j,k}^{0}|+|(\mathbf{v}_y)_{i,j,k}^{0}|+ |(\mathbf{v}_z)_{i,j,k}^{0}|)}.
\end{equation}
In practice, we take the upper bound and set the timestep as
\begin{equation}
\delta t = \frac{2h}{3 \max_{i,j,k}(|(\mathbf{v}_x)_{i,j,k}^{0}|+|(\mathbf{v}_y)_{i,j,k}^{0}|+ |(\mathbf{v}_z)_{i,j,k}^{0}|)}.
\end{equation}

In the continuous formulation, as $t \to \infty$, the density $\rho(\mathbf{x},t)$ is equalized over $D$, and the associated reference map ${\bm \xi}_{\text{final}}(\mathbf{x}) = \bm \xi(\mathbf{x},\infty)$ becomes a density-equalizing reference map. In the discrete case, the result converges to a density-equalizing reference map as the change of the density becomes very small. Therefore, we use the following convergence criterion:
\begin{equation}
\frac{\|\rho^n - \rho^{n-1}\|_2}{\text{mean}(\rho^{n-1})} \leq \epsilon,
\end{equation}
where $\epsilon$ is the error threshold. Our algorithm is summarized in Algorithm \ref{alg:3d}.

\begin{algorithm}[h]
\KwIn{A solid domain $D \subset \mathbb{R}^3$ of size $L \times M \times N$, a prescribed density $\rho^0$, the error threshold $\epsilon$, the maximum number of iterations allowed $n_{\max}$.}
\KwOut{A density-equalizing reference map ${\bm \xi}_{\text{final}}$.}
\BlankLine

Set $\delta t = \frac{2h}{3 \max_{i,j,k}(|(\mathbf{v}_x)_{i,j,k}^{0}|+|(\mathbf{v}_y)_{i,j,k}^{0}|+ |(\mathbf{v}_z)_{i,j,k}^{0}|)}$\;

Compute $A = I - \delta t \Delta$\;

Set $n = 0$\;

\Repeat{$\frac{\|\rho^n - \rho^{n-1}\|_2}{\text{mean}(\rho^{n-1})} \leq \epsilon$ or $n \geq n_{\max}$}{ 

Update $n = n + 1$\;

Solve $\rho^{n} = A^{-1} \rho^{n-1}$ as discussed in Eq.~\eqref{eqt:derm_diffusion_discrete_simplified} with an appropriate boundary condition\;

Compute the velocity field $\mathbf{v}^n$ using Eq.~\eqref{eqt:derm_velocity_discrete}\;

Update the reference map $\mathbf{\xi}^n$ using Eq.~\eqref{eqt:derm_advection_discrete}\;

}

Obtain ${\bm \xi}_{\text{final}} = {\bm \xi}^n$\;
\caption{Volumetric density-equalizing reference map}
\label{alg:3d}
\end{algorithm}

The density-equalizing reference map provides us with the information of the diffusion-based deformation. In particular, if we denote ${\bm \xi}_{\text{final}} = ({\xi}_1, {\xi}_2, {\xi}_3)$, one can obtain the forward mapping of the grid points of $D$ by tracking the intersections of the contour planes ${\xi}_1 = ih$, ${\xi}_2 = jh$ and ${\xi}_3 = kh$. More specifically, after computing the density-equalizing reference map, we build an interpolant for each of the $x$-, $y$-, and $z$-coordinates of the undeformed grid using $({\xi}_1, {\xi}_2, {\xi}_3)$, where the $\xi$ values are the abscissa and the $x,y,z$ values are the ordinate. In other words, we have three interpolant functions $I_1(x,y,z),  I_2(x,y,z), I_3(x,y,z)$ where
\begin{equation}
\left\{
\begin{split}
I_1({\xi}_1(ih,jh,kh), {\xi}_2(ih,jh,kh), {\xi}_3(ih,jh,kh)) &= ih,\\
I_2({\xi}_1(ih,jh,kh), {\xi}_2(ih,jh,kh), {\xi}_3(ih,jh,kh)) &= jh,\\
I_3({\xi}_1(ih,jh,kh), {\xi}_2(ih,jh,kh), {\xi}_3(ih,jh,kh)) &= kh,
\end{split} \right.
\end{equation} 
for all $i,j,k$. The three interpolants allow us to compute a forward deformation of the given domain. In practice, we compute a Delaunay triangulation of the given points and then perform a linear interpolation in each tetrahedral element. It is noteworthy that throughout the density-equalizing iterations it suffices to work on a grid, and the interpolation is only needed once at the end. On the contrary, a direct extension of the original density-equalizing map~\cite{Gastner04} would involve interpolating the velocity field at every iteration. Therefore, the use of the reference map technique here is advantageous.

\subsection{Boundary conditions}
Note that boundary conditions are needed for solving the diffusion equation~\eqref{eqt:derm_diffusion_discrete}. Below, we discuss three possible boundary conditions.

\begin{figure}[t]
\centering
\includegraphics[width=\textwidth]{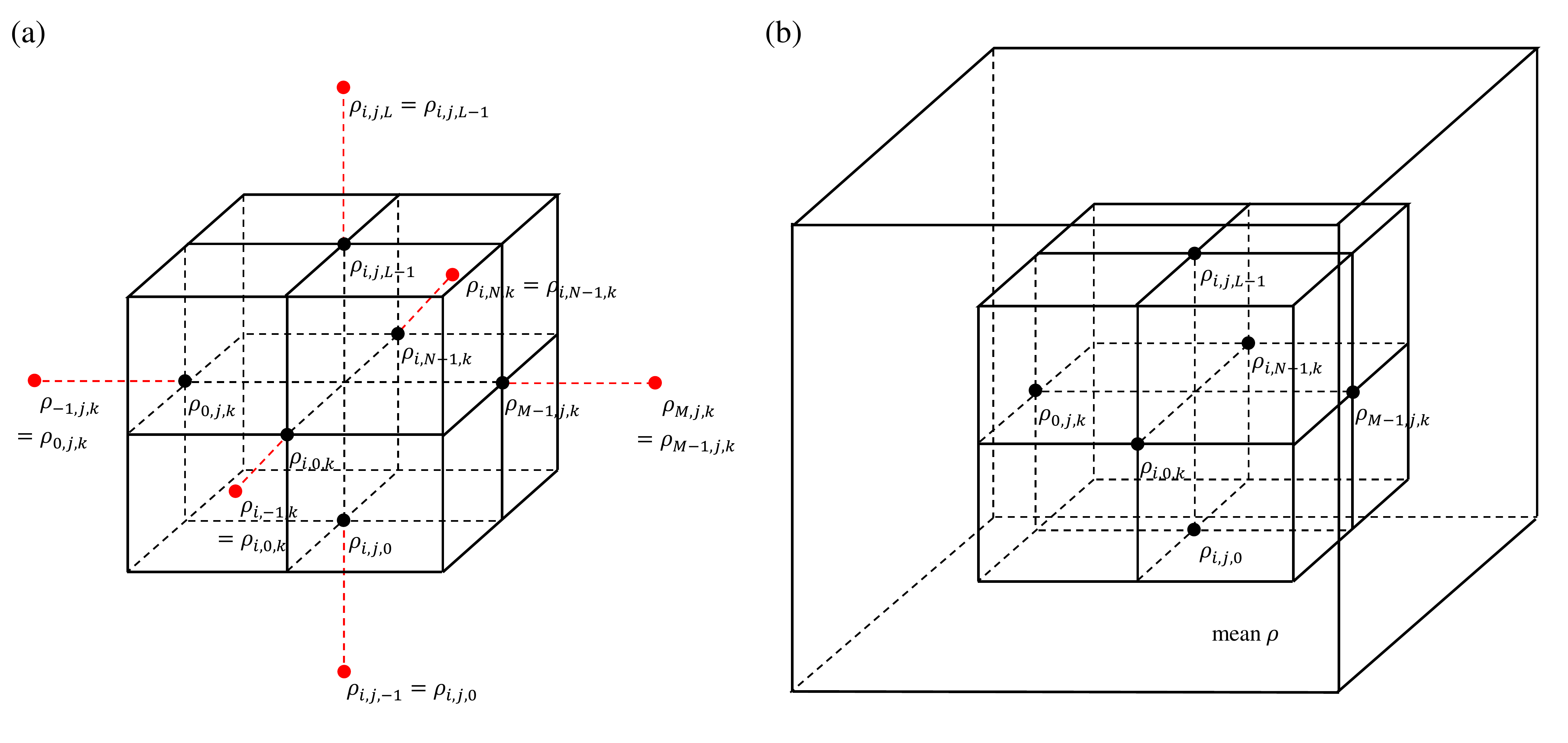}
\caption{Different choices of boundary conditions for a solid domain $D$. (a) The no-flux boundary condition can be imposed using the ghost node approach shown. A ghost node (red) is placed at the neighborhood of each boundary node on $\p D$, with the density being the same as the density at the boundary node. This ensures that there is no density flux orthogonal to the boundary faces. (b) The free boundary condition can be imposed by constructing a larger solid domain $\widetilde{D}$ and setting the density at the ``sea'' (i.e. $\widetilde{D} \setminus D$) to be the average value of the original density $\rho$. This ensures that the boundary $\partial D$ of the original domain $D$ can deform freely.}
\label{fig:derm_bc_illustration}
\end{figure}

\subsubsection{No-flux boundary condition}
To maintain the rectangular boundary shape of the 3D domain, we can enforce the no-flux boundary condition
\begin{equation}
\mathbf{n} \cdot \nabla \rho = 0
\end{equation}
at the boundaries, where $\mathbf{n}$ is the outward unit normal at the boundaries. More explicitly, this condition ensures that the density gradient at the boundaries is zero in the normal direction, and hence all the boundary planes remain planar throughout the density-equalization process.

In the discrete case, we incorporate the above no-flux boundary condition in the diffusion equation~\eqref{eqt:derm_diffusion_discrete} using the following ghost node approach. At the six boundary faces $i = 0$, $i = L-1$, $j = 0$, $j = M-1$, $k = 0$, $k = N-1$, we respectively replace the terms $\rho_{i-1,j,k}^n$, $\rho_{i+1,j,k}^n$, $\rho_{i,j-1,k}^n$, $\rho_{i,j+1,k}^n$, $\rho_{i,j,k-1}^n$, $\rho_{i,j,k+1}^n$ on the right hand side in the diffusion equation~\eqref{eqt:derm_diffusion_discrete} by $\rho_{i,j,k}^n$. This ensures that there is no density flux orthogonal to the six boundary faces throughout the density diffusion process (Fig.~\ref{fig:derm_bc_illustration}(a)).

The effect is that for each boundary node, one of its coordinates will remain unchanged under the update by the advection equation~\eqref{eqt:derm_advection_discrete}, while the other coordinates can vary. Hence, the boundary nodes can freely slide along the six boundary faces to achieve density equalization while $D$ will remain to be a rectangular box throughout the iterations.

\subsubsection{Free boundary condition}
In case it is desirable to let the domain deform freely without any boundary constraints, the ``sea'' approach by Gastner and Newman~\cite{Gastner04} can be used. Consider putting the entire solid domain $D$ with the prescribed density $\rho^0$ at the center of a larger rectangular solid domain $\widetilde{D}$. Define the density $\widetilde{\rho}^0$ on $\widetilde{D}$ by
\begin{equation}
\widetilde{\rho}^0 = \left\{\begin{array}{ll}
\rho^0 & \text{ on } D, \\
\text{mean}(\rho^0) & \text{ on } \widetilde{D} \setminus D. 
\end{array}\right.
\end{equation}
Now, we impose the no-flux boundary condition on $\p \widetilde{D}$ and compute the density-equalizing map on $\widetilde{D}$ (Fig.~\ref{fig:derm_bc_illustration}(b)). As there is no boundary condition imposed on $\p D$, $D$ can deform freely under the deformation of $\widetilde{D}$ without any components of the reference map ${\bm \xi}$ on $\p D$ being pinned. Note that the purpose of setting $\widetilde{\rho}^0 = \text{mean}(\rho^0)$ at the ``sea'' $\widetilde{D} \setminus D$ is to prevent $D$ from expanding infinitely. 

\subsubsection{Mixed boundary condition}
It is also possible to combine different boundary conditions to achieve other deformations. For instance, in case it is desired to keep the top and the bottom boundary planes planar while allowing the other boundary planes to deform freely, one can place the $L \times M \times N$ domain in a larger $\tilde{L} \times \tilde{M} \times N$ domain, where $\tilde{L} > L$ and $\tilde{M} > M$. Then, using the above-mentioned ``sea'' approach, one can enforce the no-flux boundary condition on the $z$-boundary planes and the free boundary condition on the $x$- and $y$-boundary planes, thereby producing the desired effect.

\subsection{Numerical results} 
The proposed iterative scheme is implemented in C++ with OpenMP parallelization (with grid size $h = 1$, maximum number of iterations $n_{\text{max}} = 10000$, and error threshold $\epsilon = 10^{-2}$). The sparse linear systems are solved using the conjugate gradient method (\texttt{ConjugateGradient}) in the C++ library Eigen. All experiments are performed on a PC with an Intel i7-6700K quad-core processor and 16~GB RAM. The statistics and the visualization are done using MATLAB. In the following, we assess the performance of the proposed method with numerical experiments.

\begin{figure}[t!]
\centering
\includegraphics[width=0.85\textwidth]{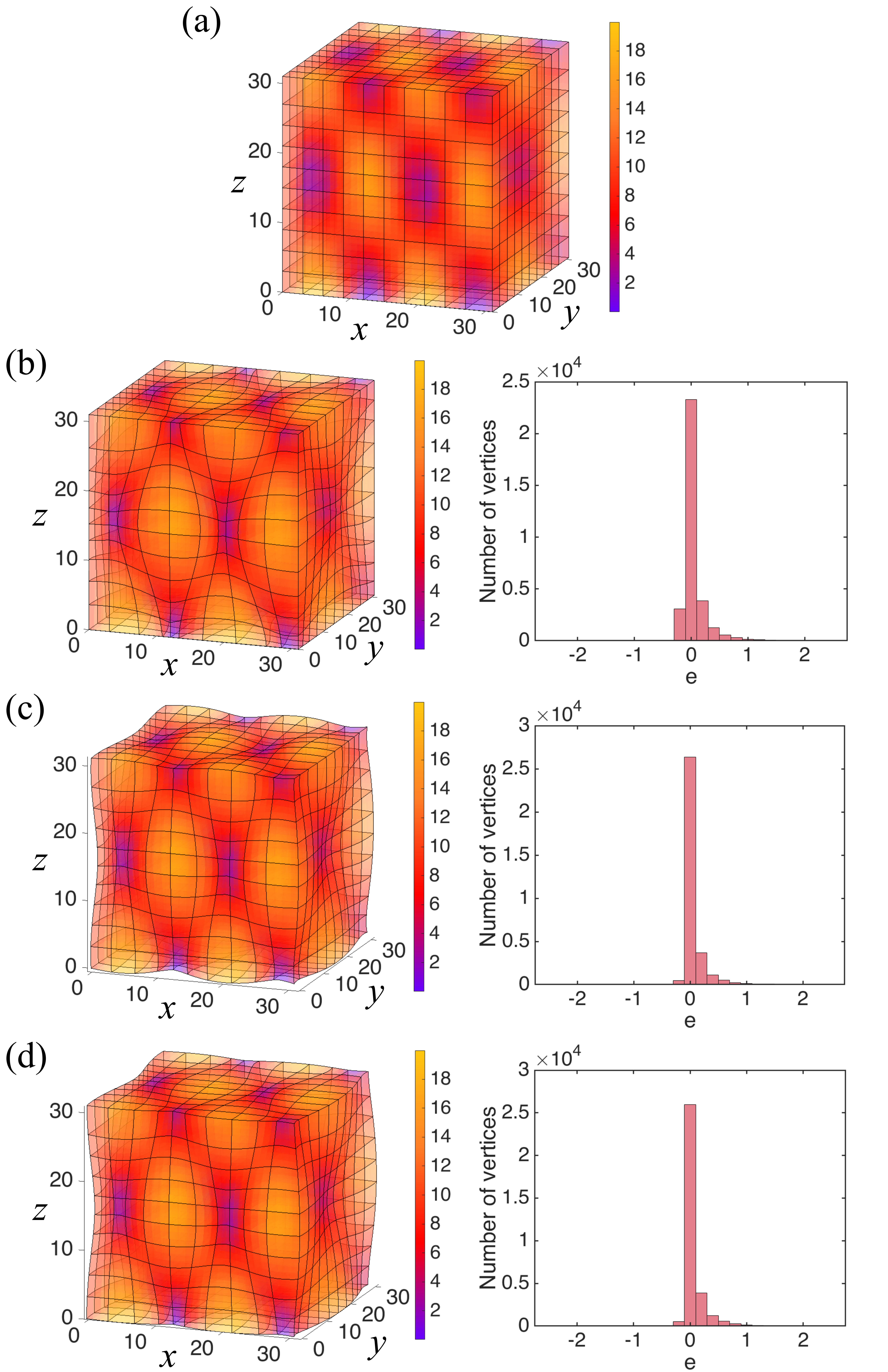}
\caption{An example of volumetric density-equalizing maps with the prescribed density being a periodic function with multiple peaks on a $32 \times 32 \times 32$ grid. The deformed contour planes at the $x$-, $y$-, and $z$-directions are plotted for visualizing the mapping results. (a) The initial state colored with the prescribed density. (b) The forward mapping result with the no-flux boundary condition and the histogram of the volume-density mismatch error $e$. (c) The forward mapping result with the free boundary condition and the histogram of $e$. (d) The forward mapping result with the mixed boundary condition where only the top and bottom boundary planes are required to remain planar, and the histogram of $e$. All the states are colored with the prescribed density. For better visualization, only part of the grid lines are drawn.}
\label{fig:3d_peaks}
\end{figure}

We first consider a cubic solid domain $D$ of grid size $L \times M \times N = 32\times 32 \times 32$ with a smooth input density 
\begin{equation}\label{eqt:peaks}
\rho^0(i,j,k) = 10+9.99\sin \left(\frac{4\pi i}{L-1}\right) \cos \left(\frac{2\pi j}{M-1}\right) \cos \left(\frac{2\pi k}{N-1}\right),
\end{equation}
with $i,j,k = 0, 1, \dots, 31$, as shown in Fig.~\ref{fig:3d_peaks}(a). We compute the volumetric density-equalizing maps with the no-flux boundary condition (Fig.~\ref{fig:3d_peaks}(b)), the free boundary condition (Fig.~\ref{fig:3d_peaks}(c)), and the mixed-boundary condition (Fig.~\ref{fig:3d_peaks}(d)). The one with the free boundary condition is done by placing the solid in a circumscribed $48 \times 48 \times 48$ cubic grid, and the one with the mixed boundary condition is achieved with the aid of a circumscribed $48 \times 48 \times 32$ rectangular grid. From Fig.~\ref{fig:3d_peaks}(b), it can be observed that different regions are enlarged or shrunk according to the prescribed density. Also, the deformed domain remains a cube, while the boundary nodes are free to slide on the boundary planes to achieve density-equalization. For the map with the free boundary condition as shown in Fig.~\ref{fig:3d_peaks}(c), it can be observed that the domain deforms freely while not maintaining a cubic shape. For the map with the mixed boundary condition as shown in Fig.~\ref{fig:3d_peaks}(d), it can be observed that the top and the bottom of the solid domain remain to be planar, while the other sides of the domain are deformed freely. 

To assess the accuracy of the mappings, we define the \emph{volume-density mismatch error} to be the logged ratio of the volumetric scale factor to the prescribed density with a normalization:
\begin{equation}
e(\mathbf{x}) = \log \frac{{\det(F(\mathbf{x}))} / \iiint_D {\det(F(\mathbf{x}))}}{ {\rho^0({\bm \xi}_{\text{final}}(\mathbf{x}))}{ / \iiint_D \rho^0}},
\end{equation}
where $F(\mathbf{x}) = \begin{pmatrix}
\frac{\partial {\bm \xi}_{\text{final}}}{\partial \mathbf{x}}
\end{pmatrix}^{-1}$ is the Jacobian, and $\iiint_D \det(F(\mathbf{x}))$ and $\iiint_D \rho^0$ are computed using the trapezoidal integration method. Note that $e = 0$ if and only if the final volume distribution matches the prescribed density distribution. As shown in Fig.~\ref{fig:3d_peaks}, the histograms of $e$ for all different boundary conditions highly concentrate at $0$, implying that the computation is accurate regardless of the choice of the boundary conditions.

\begin{figure}[t!]
\centering
\includegraphics[width=0.85\textwidth]{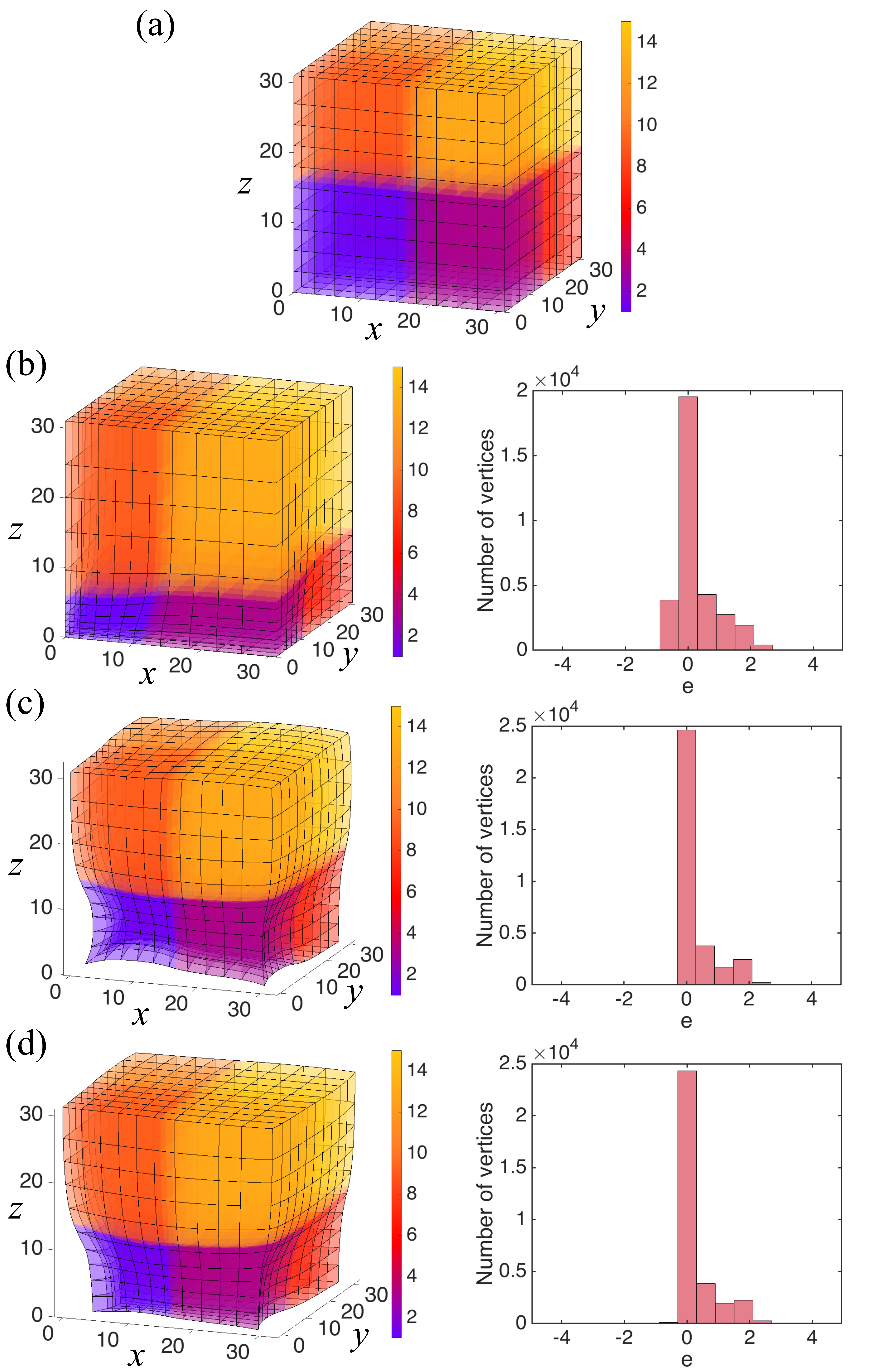}
\caption{An example of volumetric density-equalizing maps with 8 different density values defined on 8 regions on a $32 \times 32 \times 32$ grid. (a) The initial state colored with the prescribed density. (b) The forward mapping result with the no-flux boundary condition and the histogram of the volume-density mismatch error $e$. (c) The forward mapping result with the free boundary condition and the histogram of $e$. (d) The forward mapping result with the mixed boundary condition where only the top and bottom boundary planes are required to remain planar, and the histogram of $e$. All the states are colored with the prescribed density. For better visualization, only part of the grid lines are drawn.}
\label{fig:3d_8regions}
\end{figure}

Our method can also handle prescribed densities with discontinuity. Fig.~\ref{fig:3d_8regions}(a) shows another example of a 3D grid of size $L \times M \times N = 32 \times 32 \times 32$ with eight different density values defined on eight regions
\begin{equation}\label{eqt:3d_8regions}
\rho^0(i,j,k) = \left\{ \begin{array}{ll}
1 & \text{ if } i < \nicefrac{L}{2}, j < \nicefrac{M}{2}, k < \nicefrac{N}{2},\\
3 & \text{ if } i \geq \nicefrac{L}{2}, j < \nicefrac{M}{2}, k < \nicefrac{N}{2},\\
5 & \text{ if } i < \nicefrac{L}{2}, j \geq \nicefrac{M}{2}, k < \nicefrac{N}{2},\\
7 & \text{ if } i \geq \nicefrac{L}{2}, j \geq \nicefrac{M}{2}, k < \nicefrac{N}{2},\\
9 & \text{ if } i < \nicefrac{L}{2}, j < \nicefrac{M}{2}, k \geq \nicefrac{N}{2},\\
11 & \text{ if } i \geq \nicefrac{L}{2}, j < \nicefrac{M}{2}, k \geq \nicefrac{N}{2},\\
13 & \text{ if } i < \nicefrac{L}{2}, j \geq \nicefrac{M}{2}, k \geq \nicefrac{N}{2},\\
15 & \text{ if } i \geq \nicefrac{L}{2}, j \geq \nicefrac{M}{2}, k \geq \nicefrac{N}{2}.
\end{array}\right.
\end{equation}
Again, we consider computing the volumetric density-equalizing maps with the no-flux boundary condition (Fig.~\ref{fig:3d_8regions}(b)), the free boundary condition (Fig.~\ref{fig:3d_8regions}(c)) and the mixed boundary condition (Fig.~\ref{fig:3d_8regions}(d)). The one with the free boundary condition is done by placing the solid in a larger $48 \times 48 \times 48$ grid, and the one with the mixed boundary condition is achieved with the aid of a $48 \times 48 \times 32$ grid. In all cases, it can be observed that the regions are appropriately magnified or shrunk according to the input density, either with the cubic boundary shape preserved (Fig.~\ref{fig:3d_8regions}(b)), with the boundary shape freely deformed (Fig.~\ref{fig:3d_8regions}(c)), or with the planarity of only the top and the bottom boundaries preserved (Fig.~\ref{fig:3d_8regions}(d)). Also, while the input density is highly discontinuous, it can be observed from the histograms of the volume-density mismatch error $e$ that the error is concentrated at $0$. This suggests that our method is capable of handling discontinuous input densities.

\begin{table}[t]
\centering
\begin{tabular}{|c|c|c|c|} \hline
 Number of vertices ($L \times L \times L$) & Time (s) & \# iterations & $\text{mean}(|e|)$ \\ \hline
  $8 \times 8 \times 8 = 512$ & 0.02 & 105 & 0.8868 \\ \hline
  $16 \times 16 \times 16 = 4096$ & 0.20 & 323 & 0.4468 \\ \hline
  $24 \times 24 \times 24 = 13824$ & 1.14 & 504 & 0.2774 \\ \hline
  $32 \times 32 \times 32 = 32768$ & 4.43 & 637 & 0.1972 \\ \hline
  $40 \times 40 \times 40 = 64000$ & 12.13 & 727 & 0.1577 \\ \hline
  $48 \times 48 \times 48 = 110592$ & 26.65 & 784 & 0.1356 \\ \hline
  $56 \times 56 \times 56 = 175616$ & 59.59 & 816 & 0.1262 \\ \hline
  $64 \times 64 \times 64 = 262144$ & 111.53 & 856 & 0.1241 \\ \hline
\end{tabular}
\caption{The performance of the proposed method in terms of the time required for the iterative scheme, the number of iterations needed, and the mean of the absolute volume-density mismatch error $\text{mean}(|e|)$. In all experiments, four OpenMP threads are used.}
\label{table:performance}
\end{table}

We further examine the performance of the proposed method with different resolutions $L \times L \times L$ with the density given by Eq.~\eqref{eqt:peaks}. For a fair assessment, we introduce a diffusion coefficient $\kappa = L/64$ on the right-hand side of Eq.~\eqref{eqt:derm_diffusion} and Eq.~\eqref{eqt:velocity_field} to ensure that the rate of diffusion is proportionally rescaled. From Table~\ref{table:performance}, it can be observed that the volume-density mismatch error $\text{mean}(|e|)$ is inversely proportional to $L$. Therefore, the accuracy of the density-equalizing reference map can be improved by increasing the resolution of the numerical domain. 

\section{Applications} \label{sect:application}
\subsection{Volumetric data visualization}
Recall that the original density-equalizing map for 2D domains was developed primarily for sociological and biological data visualization on the world map. Instead of simply representing certain information using colors on a normal world map, one can convey the information geometrically by distorting the map via density equalization. Analogously, with the aid of our proposed VDERM method, we can provide an alternative, geometrical visualization for volumetric data. We demonstrate this idea using two medical and sociological examples. 

\subsubsection{Medical data visualization}
\begin{figure}[t!]
\centering
\includegraphics[width=0.6\textwidth]{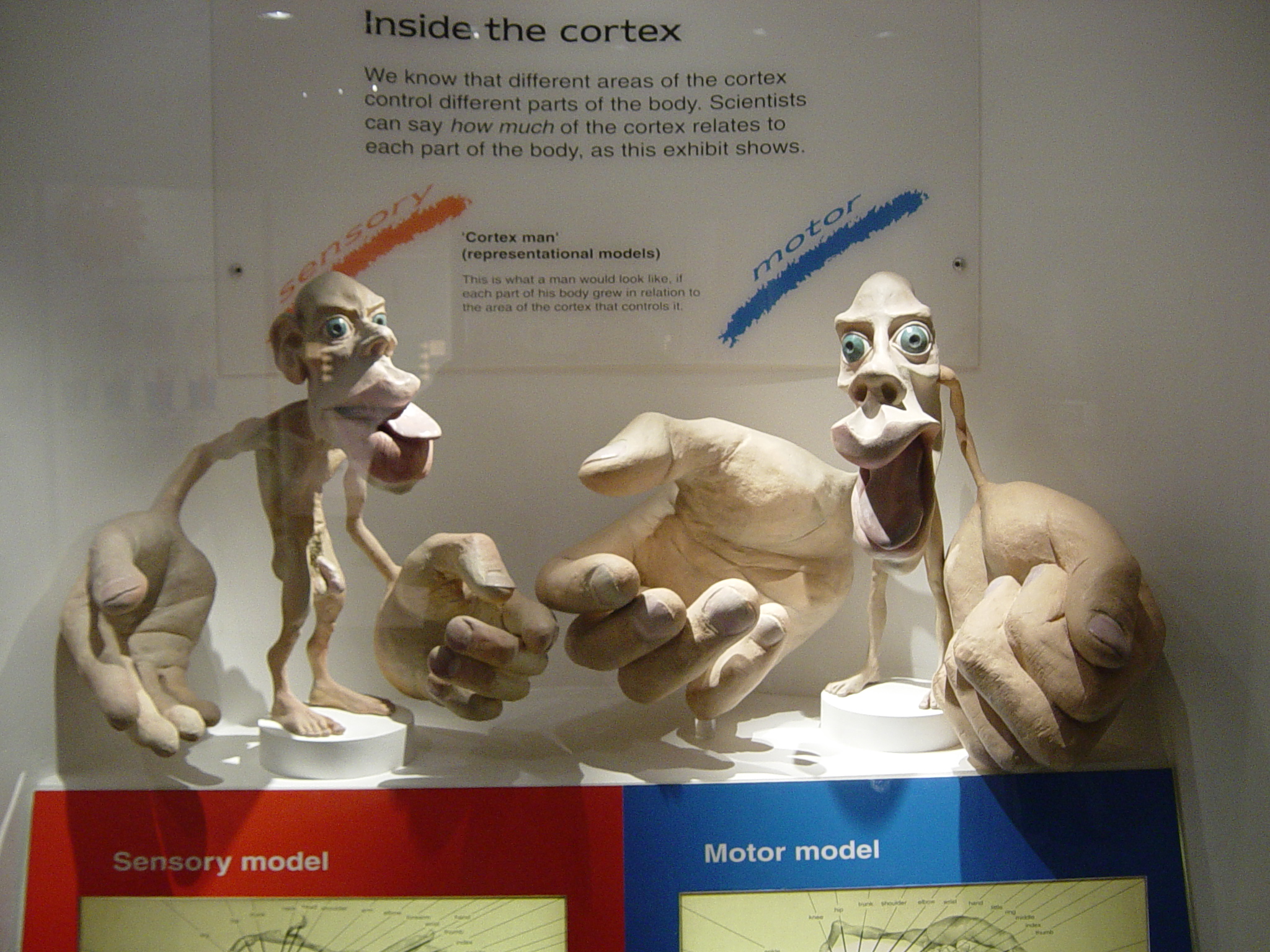}
\caption{Sculptures of sensory and motor homunculi at the Museum of Natural History, London, UK. The image is adopted online~\cite{sensory} under the CC BY-SA 3.0 license.}
\label{fig:homunculi_physical}
\end{figure}

In neurology, a \emph{cortical homunculus} (also known as a \emph{cortex man}) is a 3D human model with distortions at different parts of the body representing the proportion of the human brain that is used for processing sensory or motor functions for those parts~\cite{Marieb07human}. Fig.~\ref{fig:homunculi_physical} shows sculptures of sensory and motor homunculi at the Museum of Natural History, London, UK. Such distorted representations enable us to easily understand the ratios between the levels of nerve control for different parts of the human body.

Our proposed VDERM method is well-suited for this medical data visualization task. To illustrate this idea, we adopt medical data of the spatial acuity for pain in human body from Ref.~\cite{Mancini14}. The spatial acuity was assessed by measuring the 2-point discrimination (2PD) thresholds (in cm) at different parts of the body, where a smaller threshold implies a higher spatial acuity. To represent the data geometrically, we deform a 3D human body model using our proposed volumetric density-equalizing reference map method, with the density being the reciprocal of the 2PD threshold of each part of the body divided by the volume of it:
\begin{equation}
\rho^0 = \left\{\begin{array}{ll}
\nicefrac{(1/0.6)}{\text{Volume of the fingers}} & \text{ for the fingers,}\\
\nicefrac{(1/1.1)}{\text{Volume of the hand palms}} & \text{ for the palms,}\\
\nicefrac{(1/2.7)}{\text{Volume of the dorsum of the hands}} & \text{ for the dorsum of the hands,}\\
\nicefrac{(1/2.1)}{\text{Volume of the arms}} & \text{ for the arms,}\\
\nicefrac{(1/1.2)}{\text{Volume of the head}} & \text{ for the head,}\\
\nicefrac{(1/1.4)}{\text{Volume of the shoulders}} & \text{ for the shoulders,}\\
\nicefrac{(1/2.5)}{\text{Volume of the lower back}} & \text{ for the lower back,}\\
\nicefrac{(1/2.9)}{\text{Volume of the thighs}} & \text{ for the thighs,}\\
\nicefrac{(1/3.3)}{\text{Volume of the calves}} & \text{ for the calves,}\\
\nicefrac{(1/3.5)}{\text{Volume of the dorsum of the foot}} & \text{ for the dorsum of the foot,}\\
\nicefrac{(1/1.3)}{\text{Volume of the foot soles}} & \text{ for the foot soles.}\\
\end{array}\right.
\end{equation}
The density at the other parts of the body is defined using the closest parts with data available. The density at the remaining regions of the volumetric domain is set to be the average density of the human body, making those regions a ``sea'' to allow for a free-boundary deformation of the human body model.

\begin{figure}[t!]
\centering
\includegraphics[width=0.8\textwidth]{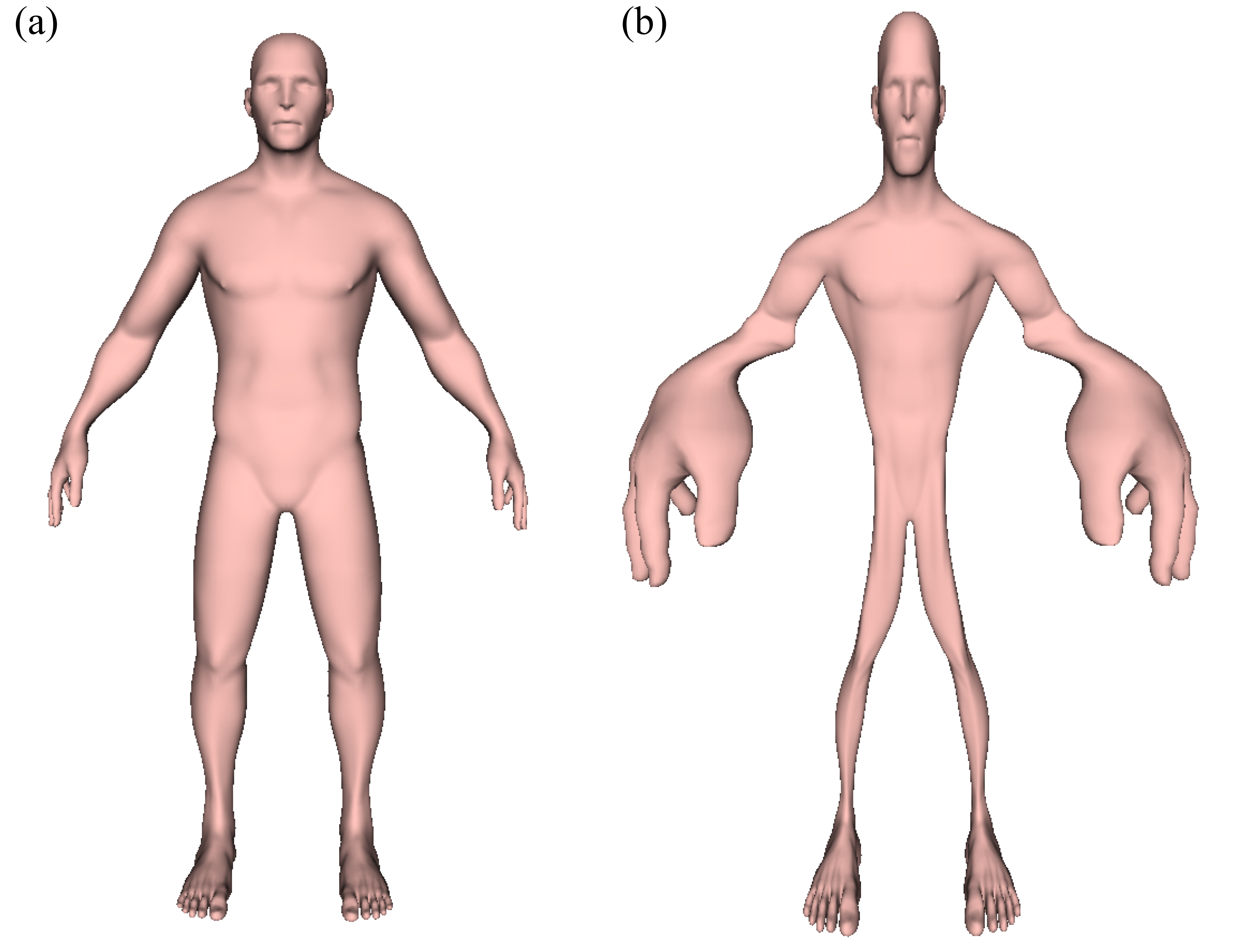}
\caption{Producing a cortical homunculus using our proposed VDERM method. (a) A normal 3D human body model. (b) A deformed human model produced by our proposed method representing the spatial acuity of pain~\cite{Mancini14}. }
\label{fig:homunculi_numerical}
\end{figure}

Fig.~\ref{fig:homunculi_numerical} shows the deformed human body model obtained by our proposed VDERM method. It can be observed that the hands are expanded significantly, as they are with the highest spatial acuity. The deformed model and the sculptures in Fig.~\ref{fig:homunculi_physical} are highly similar, except for the mouth/lips at which the acuity data is not provided in Ref.~\cite{Mancini14}. This example shows that our proposed method is useful for medical data visualization. 

\begin{figure}[t!]
\centering
\includegraphics[width=0.65\textwidth]{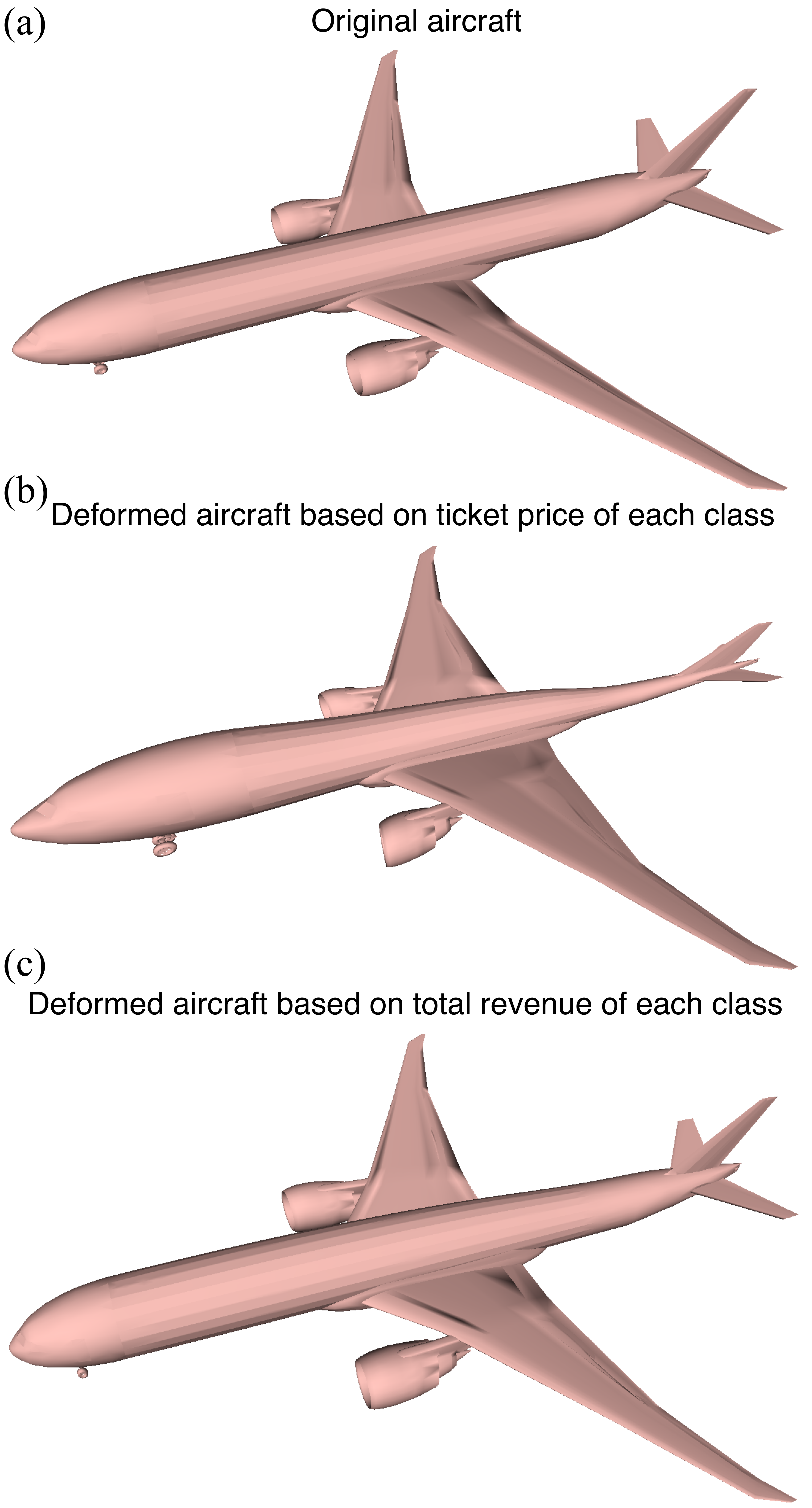}
\caption{Visualization of the economics of airline class using our proposed VDERM method. (a) A 3D model of the Boeing 777-300ER aircraft. (b) The deformed aircraft produced by our proposed method based on the ticket price for each ticket class. The magnification at the front and the shrinkage at the end (in $x$-, $y$-, and $z$-directions) reflect the large ticket price difference between the first class seats and the economy class seats. (c) The deformed aircraft produced by our proposed method based on the total revenue for each ticket class. The relatively small deformation suggests that the total revenue of each ticket class, with the amount of cabin space taken into account, is similar for the four classes.}
\label{fig:airline}
\end{figure}
\subsubsection{Sociological data visualization}
Our proposed VDERM method can also be used for visualizing sociological data. For instance, airlines have multiple flight classes, where passengers of different flight classes enjoy service and accommodation  at different levels. It is well known that the ticket price of different travel classes can be significantly different. Here, we apply our proposed method for visualizing the ticket price of different travel classes on an aircraft. We collect the ticket price data of a round-trip direct flight between New York (JFK) and Hong Kong (HKG) (departure date: March 1, 2020; returning date: March 8, 2020; data retrieved on December 3, 2019 from the American Airlines website~\cite{aa}). Both the departure flight and the returning flight are operated on a Boeing 777-300ER aircraft, in which there are in total 6 first class seats, 53 business class seats, 34 premium economy class seats and 182 economy class seats~\cite{seatguru}. The ticket prices for a first, business, premium economy, and economy class seat are \$16039, \$6922, \$2542, and \$941 respectively. 

To give a geometrical interpretation of the ticket price difference, we deform a 3D model of the Boeing 777-300ER aircraft using our proposed volumetric density-equalizing reference map method, with the density being the ticket price divided by the volume of the cabin:
\begin{equation}
\rho^0 = \left\{\begin{array}{ll}
\nicefrac{16039}{\text{Volume of the first class cabin}} & \text{ for the first class cabin,}\\
\nicefrac{6922}{\text{Volume of the business class cabin}} & \text{ for the business class cabin,}\\
\nicefrac{2542}{\text{Volume of the premium economy class cabin}} & \text{ for the premium economy class cabin,}\\
\nicefrac{941}{\text{Volume of the economy class cabin}} & \text{ for the economy class cabin.}
\end{array}\right.
\end{equation}
Again, the density at the remaining regions of the volumetric domain is set to be the average density of the above to allow for a free-boundary deformation. Note that the division by the volume of the cabins is to ensure that the ticket price ratio is equal to the ratio of the volume of the entire cabins. Fig.~\ref{fig:airline}(a),(b) show the original Boeing 777-300ER aircraft and the deformed aircraft. It can be observed that the first class cabin expands in all directions, while the economy class cabin shrinks significantly. This reflects the large ticket price difference between the first class seats and the economy class seats. 

We can also estimate the total revenue of different ticket classes by multiplying the ticket price by the total number of seats in each class. This gives a total revenue of \$96234, \$366866, \$86428, and \$171262 for the four classes respectively. Similarly, we obtain a volumetric density-equalizing map using the following density:
\begin{equation}
\rho^0 = \left\{\begin{array}{ll}
\nicefrac{96234}{\text{Volume of the first class cabin}} & \text{ for the first class cabin,}\\
\nicefrac{366866}{\text{Volume of the business class cabin}} & \text{ for the business class cabin,}\\
\nicefrac{86428}{\text{Volume of the premium economy class cabin}} & \text{ for the premium economy class cabin,}\\
\nicefrac{171262}{\text{Volume of the economy class cabin}} & \text{ for the economy class cabin.}
\end{array}\right.
\end{equation}
Fig.~\ref{fig:airline}(c) shows the deformed aircraft. Unlike the previous result, this deformed aircraft is not significantly different from the original model. We observe that the economy class cabin is slightly shrunk, while the business class cabin is slightly expanded. This suggests that the total revenue of each ticket class, with the amount of cabin space used taken into account, is in fact similar for the four classes.

The above examples show that our proposed VDERM method is capable of producing deformations to convey sociological information in volumetric domains.

\begin{figure*}[t!]
\centering
\includegraphics[width=0.9\textwidth]{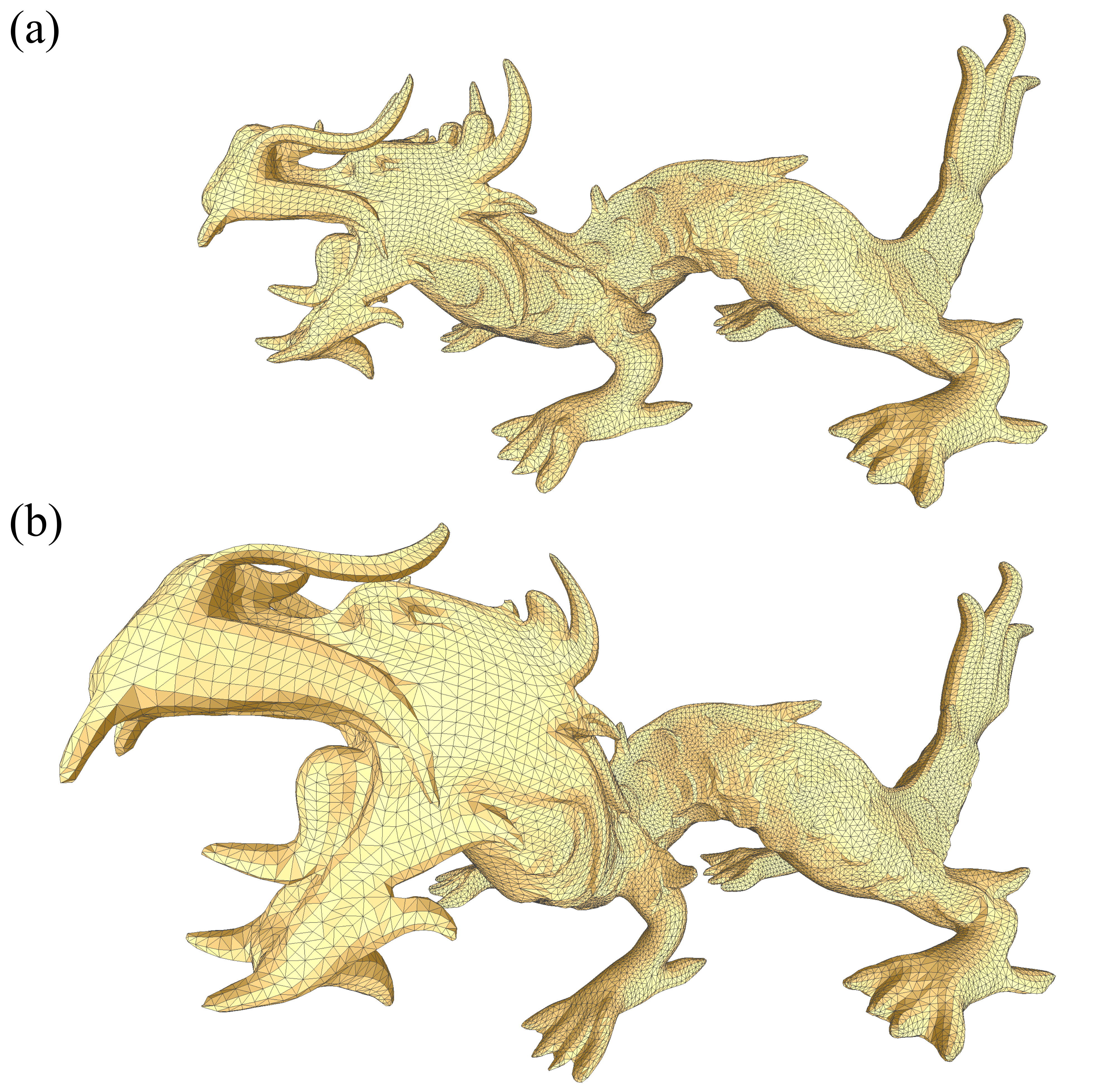}
\caption{Deformation of a dragon model using our volumetric density-equalizing reference map. (a) The original dragon. (b) The deformed dragon with the head enlarged.}
\label{fig:dragon}
\end{figure*}

\subsection{Deformation-based shape modeling}
Besides volumetric data visualization, our proposed VDERM method can also be used for deformation-based shape modeling. Suppose we are given a mesh in $\mathbb{R}^3$. By defining different density values at different parts of an underlying 3D grid and computing the volumetric density-equalizing map, we obtain a deformation of the grid, which induces a deformation of the object. Using this idea, we can easily achieve different shape modeling effects.

Fig.~\ref{fig:dragon} shows an experiment on a dragon model adapted from The Stanford 3D Scanning Repository~\cite{Stanford}. In this experiment, we aim to magnify the head of the dragon while keeping the body shape unchanged. On an underlying 3D grid of size $32\times 32\times 32$, we set the density to be
\begin{equation}
\rho^0 = \left\{\begin{array}{ll}
10 & \text{ around the head of the dragon,}\\
1 & \text{ elsewhere.}
\end{array}\right.
\end{equation}
We then compute the volumetric density-equalizing map on the underlying grid and obtain the induced deformation of the dragon. It can be observed that the head of dragon is effectively enlarged, while the main body of the dragon is almost unchanged. This demonstrates the advantage of our approach in accurately producing a deformed model that satisfies the desired effect. Also, while the input density $\rho$ is with a sharp discontinuity, the density diffusion process produces a smooth transition from the high-density region to the low-density region. As shown in Fig.~\ref{fig:dragon}, the neck of the dragon is also slightly enlarged under the deformation because of the density diffusion process. This demonstrates another advantage of our approach in producing a naturally and smoothly deformed model.

\begin{figure}[t!]
\centering
\includegraphics[width=0.75\textwidth]{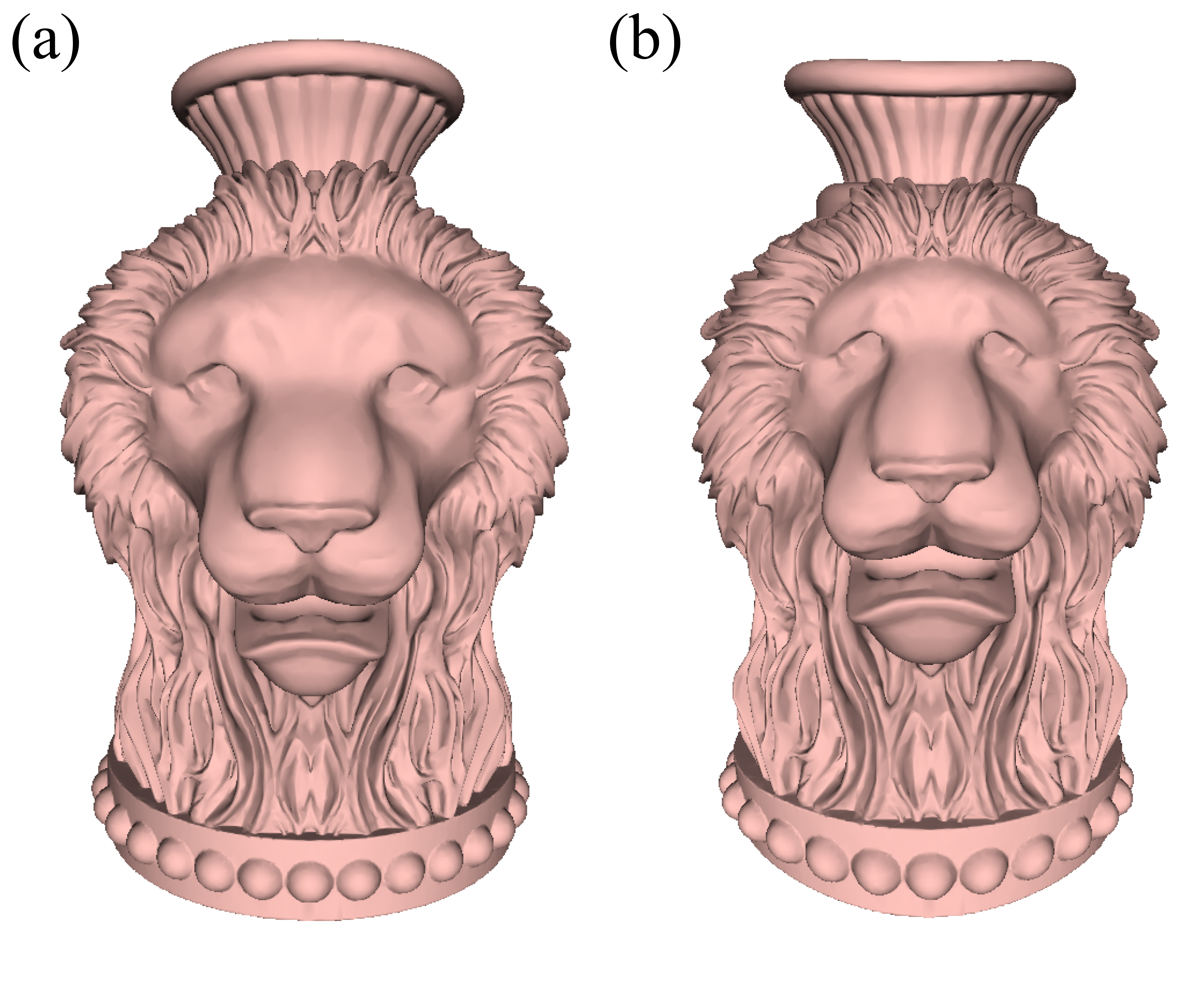}
\caption{Deformation of a lion vase model using our proposed VDERM method. (a) The original lion vase. (b)~The deformed lion vase with the facial expression changed.}
\label{fig:vaselion}
\end{figure}

Fig.~\ref{fig:vaselion} shows another experiment of shape modeling using our approach. This time, we consider a lion vase model adapted from the 3D Segmentation Benchmark~\cite{3dsegbenchmark}. To change the facial expression of the lion, we define a smaller density value around the forehead of the lion on the underlying grid and a larger density value around the chin of it. Then we compute the volumetric density-equalizing map and obtain the induced deformation of the lion vase. It can be observed that the facial expression of the lion is changed naturally under our deformation-based approach. 
 
\subsection{Shape morphing}
\begin{figure*}[t!]
\centering
\includegraphics[width=\textwidth]{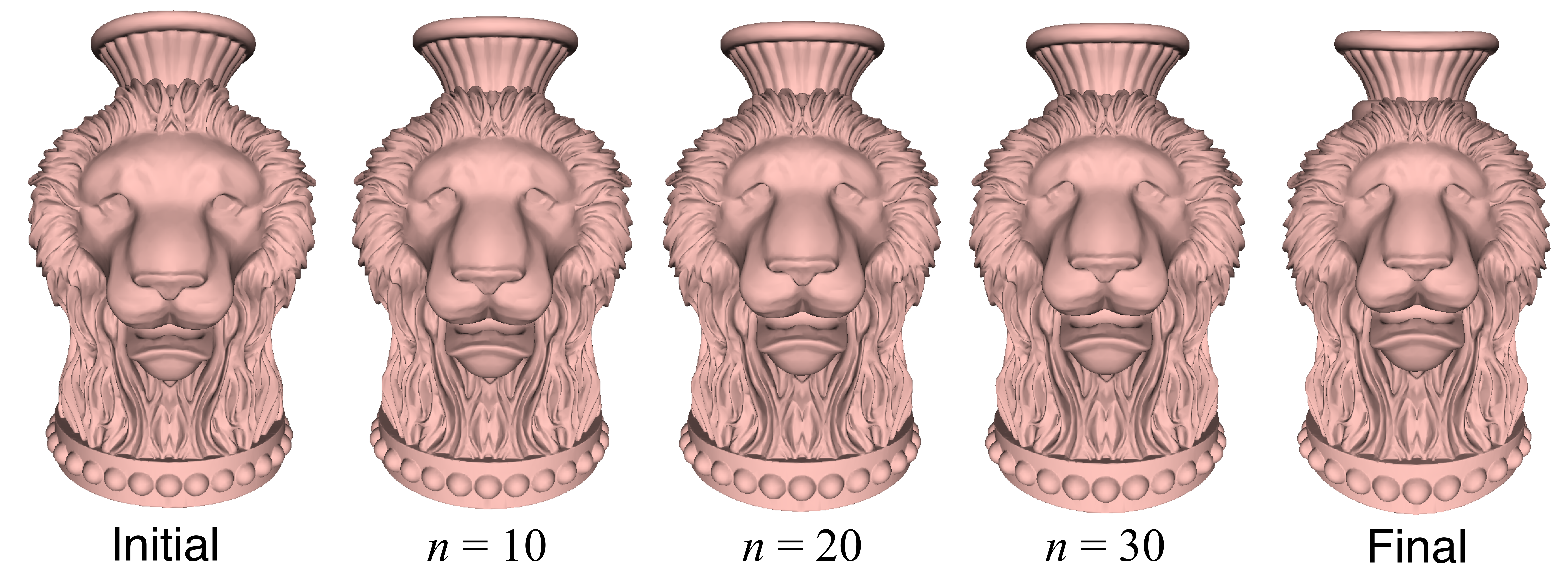}
\caption{An animation of the facial expression change of the lion produced using the intermediate states of the volumetric density-equalizing reference maps. }
\label{fig:vaselion_animation}
\end{figure*}

It is noteworthy that our proposed approach deforms the underlying grid continuously with an iterative scheme, and hence the intermediate states of the deformation of the grid can be used for producing a continuous change of the 3D model. Fig.~\ref{fig:vaselion_animation} shows multiple snapshots of the facial expression change produced by computing the induced deformation of the lion vase with respect to the intermediate states of the deformed 3D grid at different time points. It can be observed that the facial expression changes continuously. This experiment suggests that our deformation-based method is also advantageous for applications which focus not only on the final shape but also the intermediate states, such as object morphing. 

\begin{figure}[t!]
\centering
\includegraphics[width=0.85\textwidth]{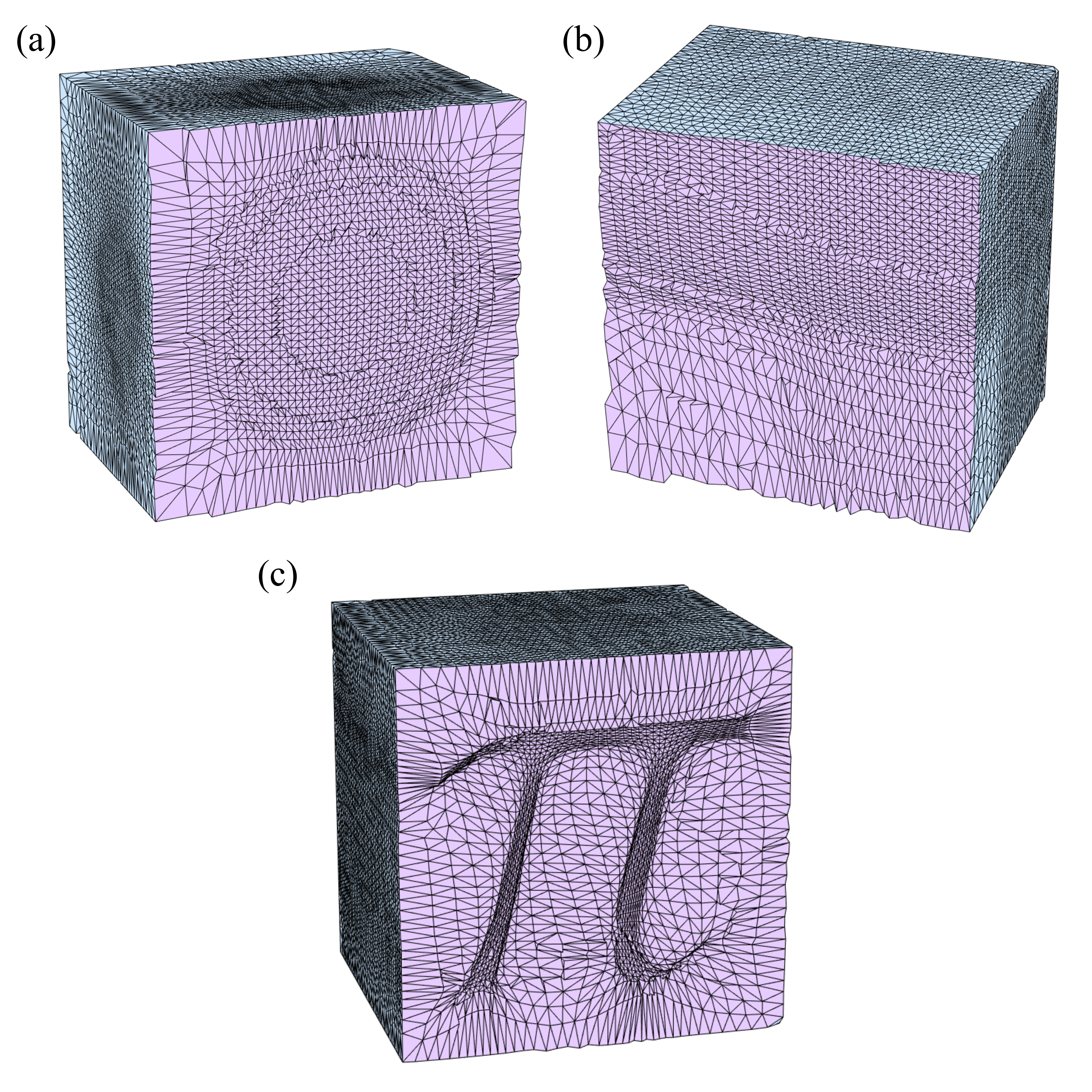}
\caption{Examples of adaptive remeshing achieved by our proposed VDERM method. (a) A cube with the central part refined. (b) A cube with eight different tetrahedral mesh densities. (c) A cube with a $\pi$ shape having a higher mesh density. For better visualization, the surface is colored in blue and the interior tetrahedral elements are colored in purple.}
\label{fig:remeshing}
\end{figure}

\subsection{Adaptive remeshing}
Our proposed VDERM method can also be applied for remeshing, which aims at improving the discretization of meshes~\cite{Choi2016spherical,Choi2016fast}. More specifically, suppose our goal is to construct a tetrahedral mesh for a genus-0 closed surface. Furthermore, the tetrahedral mesh density of some prescribed regions in the volumetric domain $D$ is desired to be higher. To achieve this, we run our proposed algorithm with a larger $\rho$ at those regions and obtain the volumetric density-equalizing reference map ${\bm \xi}_{\text{final}}$. Then, we construct a uniform tetrahedral mesh in the deformed domain using DistMesh~\cite{Persson2004a}, and map it back to the original domain $D$ using ${\bm \xi}_{\text{final}}$. This puts more tetrahedral elements in the regions with a larger prescribed $\rho$ and less elements in the regions with a smaller $\rho$, thereby achieving a tetrahedral mesh with adaptive mesh density. 

Three examples are given in Fig.~\ref{fig:remeshing}. For Fig.~\ref{fig:remeshing}(a), the central part of the domain is desired to be with a higher tetrahedral mesh density. To achieve this, we set
\begin{equation}
\rho^0 = \left\{\begin{array}{ll}
10 & \text{ around the center of the domain,}\\
1 & \text{ elsewhere.}
\end{array}\right.
\end{equation}
For Fig.~\ref{fig:remeshing}(b), we set $\rho^0$ by Eq.~\eqref{eqt:3d_8regions} in order to achieve eight different tetrahedral mesh densities at the eight corners. For Fig.~\ref{fig:remeshing}(c), we set $\rho^0$ to be 10 inside the $\pi$ shape and 1 otherwise, in order to achieve a higher mesh density at the $\pi$ shape. In all examples, the resulting volumetric density-equalizing reference map ${\bm \xi}_{\text{final}}$ is capable of mapping a uniform tetrahedral mesh to the original domain and produces the desired tetrahedral meshes with adaptive mesh density. 

DistMesh~\cite{Persson2004a} also allows for adaptive remeshing, by specifying an edge length function expression that yields the target adaptive resolution. This differs from our approach, where we can exactly specify the volumes of mesh elements throughout the domain. Our approach allows us to deform a mesh to precisely control element volume while maintaining mesh topology, although it can result in elongated triangles such as those shown in Fig.~\ref{fig:remeshing}(c). A hybrid approach combining both DistMesh and VDERM provides a flexible platform for generating meshes for numerical computations (e.g.~solving PDEs) with added precision in certain regions.

\section{Discussion} \label{sect:discussion}
In this work, we have proposed a method for computing volumetric density-equalizing reference maps with a prescribed density defined in a solid domain in $\mathbb{R}^3$. This is the first work on the formulation of density-equalizing map in 3D for diffusion-based deformation of volumetric domains. We have demonstrated the effectiveness of the proposed VDERM method via numerical experiments and discussed potential applications of the proposed method to volumetric data visualization, remeshing and deformation-based shape modeling. We have further introduced a novel use of the entire density-diffusion process for time-dependent applications such as object morphing.

One possible future work is to extend our VDERM method for achieving more sophisticated shape modeling effects. Note that our approach currently focuses on deformations related to relative volume change of voxels. We may consider adding some extra steps with other deformation energies~\cite{Hildebrandt11} throughout the density-equalization iterations to produce other type of deformations. For instance, we may enable rotations of the velocity vectors at certain regions throughout the iterations in order to achieve a rotational effect on the deformed object. It will also be interesting to explore if our method can be combined with landmark-matching mapping methods~\cite{Choi15,Lee16} to achieve a wider range of shape deformations that satisfy prescribed landmark correspondences in 3D. 

\section*{Conflict of interest}
The authors declare that they have no conflict of interest.

\section*{Acknowledgments} 
This work was supported in part by the Croucher Foundation (to G. P. T. Choi), the Harvard Quantitative Biology Initiative and the NSF-Simons Center for Mathematical and Statistical Analysis of Biology at Harvard, award number \#1764269 (to G. P. T. Choi), and the Applied Mathematics Program of the U.S. Department of Energy (DOE) Office of Science Advanced Scientific Computing Research under contract DE-AC02-05CH11231 (to C. H. Rycroft). We thank Dr.~Anselm Hui (Prince of Wales Hospital) for useful discussion on medical data visualization.


\begin{thebibliography}{}
\bibitem{Gastner04}
M. T. Gastner and M. E. J. Newman. \textit{Diffusion-based method for producing density-equalizing maps}. Proceedings of the National Academy of Sciences of the United States of America, 101(20), 7499--7504, 2004.

\bibitem{Wake08}
D. B. Wake and V. T. Vredenburg, \textit{Are we in the midst of the sixth mass extinction? A view from the world of amphibians}. Proceedings of the National Academy of Sciences of the United States of America, 105(1), 11466--11473, 2008.

\bibitem{Dorling10}
D. Dorling, M. Newman, and A. Barford, \textit{The atlas of the real world: mapping the way we live}. Thames \& Hudson, 2010.

\bibitem{Pan12}
R. K. Pan, K. Kaski, and S. Fortunato, \textit{World citation and collaboration networks: uncovering the role of geography in science}. Scientific Reports, 2, 902, 2012.

\bibitem{Matthews14}
H. D. Matthews, T. L. Graham, S. Keverian, C. Lamontagne, D. Seto, and T. J. Smith, \textit{National contributions to observed global warming}. Environmental Research Letters, 9(1), 014010, 2014.

\bibitem{Dodd16}
P. J. Dodd, C. Sismanidis, and J. A. Seddon, \textit{Global burden of drug-resistant tuberculosis in children: a mathematical modelling study}. The Lancet infectious diseases, 16(10), 1193--1201, 2016.

\bibitem{Ballas17}
D. Ballas, D. Dorling, and B. Hennig, \textit{Analysing the regional geography of poverty, austerity and inequality in Europe: a human cartographic perspective}. Regional Studies, 51(1), 174--185, 2017.

\bibitem{Choi18}
G. P. T. Choi and C. H. Rycroft, \textit{Density-equalizing maps for simply connected open surfaces}. SIAM Journal on Imaging Sciences, 11(2), 1134--1178, 2018.

\bibitem{Gastner18}
M. T. Gastner, V. Seguy, and P. More, \textit{Fast flow-based algorithm for creating density-equalizing map projections}. Proceedings of the National Academy of Sciences of the United States of America, 201712674. 2018.

\bibitem{Choi20}
G. P. T. Choi, B. Chiu, and C. H. Rycroft, \textit{Area-preserving mapping of 3D carotid ultrasound images using density-equalizing reference map}. IEEE Transactions on Biomedical Engineering, 1--11, 2020.

\bibitem{Kamrin12}
K. Kamrin, C. H. Rycroft, and J.-C. Nave, \textit{Reference map technique for finite-strain elasticity and fluid–solid interaction}. Journal of the Mechanics and Physics of Solids, 60(11), 1952--1969, 2012.

\bibitem{Valkov15}
B. Valkov, C. H. Rycroft, and K. Kamrin, \textit{Eulerian method for multiphase interactions of soft solid bodies in fluids}. Journal of Applied Mechanics, 82(4), 041011, 2015. 

\bibitem{Rycroft18}
C. H. Rycroft, C.-H. Wu, Y. Yu, and K. Kamrin, \textit{Reference Map Technique for Incompressible Fluid-Structure Interaction}. Preprint, arXiv:1810.03015.

\bibitem{Gurtin10}
M. E. Gurtin, E. Fried, and L. Anand, \textit{The mechanics and thermodynamics of continua}. Cambridge University Press, 2010.

\bibitem{Govindjee96}
S. Govindjee and P. A. Mihalic, \textit{Computational methods for inverse finite elastostatics}. Computer Methods in Applied Mechanics and Engineering, 136(1--2), 47--57, 1996.

\bibitem{Fachinotti08}
V. D. Fachinotti, A. Cardona, and P. Jetteur, \textit{Finite element modelling of inverse design problems in large deformations anisotropic hyperelasticity}. International Journal for Numerical Methods in Engineering, 74(6), 894--910, 2008.

\bibitem{Courant28}
R. Courant, K. Friedrichs, H. Lewy, \textit{\"Uber die partiellen Differenzengleichungen der mathematischen Physik}. Mathematische Annalen (in German), 100(1), 32--74, 1928.

\bibitem{Stanford}
The Stanford 3D Scanning Repository. \url{http://graphics.stanford.edu/data/3Dscanrep/}


\bibitem{Marieb07human}
E. Marieb and K. Hoehn, \textit{Human Anatomy and Physiology}. Pearson Benjamin Cummings: San Francisco, 2007.

\bibitem{Mancini14}
F. Mancini, A. Bauleo, J. Cole, F. Lui, C. A. Porro, P. Haggard, and G. D. Iannetti, \textit{Whole-body mapping of spatial acuity for pain and touch}. Annals of Neurology, 75(6), 917--924, 2014.

\bibitem{sensory}
\textit{Cortical homunculus}. \url{https://en.wikipedia.org/wiki/Cortical\_homunculus}


\bibitem{aa}
\textit{American Airlines}. \url{https://www.aa.com} 

\bibitem{seatguru}
\textit{SeatGuru}. \url{https://www.seatguru.com} 


\bibitem{3dsegbenchmark}
3D Segmentation Benchmark. \url{http://193.48.251.101/3dsegbenchmark/bust.html}

\bibitem{Choi2016spherical}
G. P.-T. Choi, K. T. Ho, and L. M. Lui, \textit{Spherical conformal parameterization of genus-0 point clouds for meshing}. SIAM Journal on Imaging Sciences, 9, 1582--1618 (2016).

\bibitem{Choi2016fast}
G. P.-T. Choi, M. H.-Y. Man, and L. M. Lui, \textit{Fast spherical quasiconformal parameterization of genus-$0$ closed surfaces with application to adaptive remeshing}. Geometry, Imaging \& Computing, 3(1--2), 1--29, 2016.

\bibitem{Persson2004a}
P.-O. Persson and G. Strang, \textit{A Simple Mesh Generator in MATLAB}. SIAM Review, 46(2), 329--345, 2004.

\bibitem{Hildebrandt11}
K. Hildebrandt, C. Schulz, C. V. Tycowicz, and K. Polthier, \textit{Interactive surface modeling using modal analysis}. ACM Transactions on Graphics (TOG), 30(5), 119, 2011.


\bibitem{Choi15}
P. T. Choi, K. C. Lam, and L. M. Lui, \textit{FLASH: Fast landmark aligned spherical harmonic parameterization for genus-0 closed brain surfaces}. SIAM Journal on Imaging Sciences, 8(1), 67--94, 2015.

\bibitem{Lee16}
Y. T. Lee, K. C. Lam, and L. M. Lui, \textit{Landmark-matching transformation with large deformation via n-dimensional quasi-conformal maps}. Journal of Scientific Computing, 67(3), 926--954, 2016.

\end{thebibliography}
\end{document}